\title{The Quasi-Randomness of Hypergraph Cut Properties}

\author{Asaf Shapira
\thanks{School of Mathematics and School of Computer Science, Georgia Institute of Technology, Atlanta, GA 30332. Supported in part by NSF Grant DMS-0901355.} \and Raphael Yuster\thanks{Department of
Mathematics, University of Haifa, Haifa 31905, Israel. E--mail:
raphy@math.haifa.ac.il} }

\date{}

\documentclass [letterpaper,11pt]{article}
\usepackage{amsfonts}

\newtheorem{theo}{Theorem}
\newtheorem{problem}{Problem}

\newtheorem{coro}[theo]{Corollary}

\newtheorem{prop}[theo]{Proposition}
\newtheorem{lemma}{Lemma}[section]
\newtheorem{definition}[lemma]{Definition}

\newcommand{\qed}{\hspace*{\fill} \rule{7pt}{7pt}}

\newcommand{\ignore}[1]{}

\begin{document}
\maketitle

\begin{abstract}

Let $\alpha_1,\ldots,\alpha_k$ satisfy $\sum_i \alpha_i=1$ and suppose a $k$-uniform hypergraph on $n$ vertices satisfies the
following property; in any partition of its vertices into $k$ sets $A_1,\ldots,A_k$ of sizes $\alpha_1 n,\ldots,\alpha_k n$, the number
of edges intersecting $A_1,\ldots,A_k$ is the number one would expect to find in a random $k$-uniform hypergraph. Can we then infer that
$H$ is quasi-random? We show that the answer is negative if and only if $\alpha_1=\cdots= \alpha_k=1/k$. This resolves an open problem raised in 1991 by Chung and Graham [J. AMS '91].

While hypergraphs satisfying the property corresponding to $\alpha_1=\cdots= \alpha_k=1/k$ are not necessarily quasi-random, we manage to find
a characterization of the hypergraphs satisfying this property. Somewhat surprisingly, it turns out that (essentially) there is a unique non quasi-random hypergraph satisfying this property. The proofs combine probabilistic and algebraic arguments with results from the theory of association schemes.

\end{abstract}

\section{Introduction}\label{intro}

We study quasi-random hypergraphs (and graphs), that is, hypergraphs which have the properties one would expect to find
in ``truly'' random hypergraphs. We focus on $k$-uniform hypergraphs $H=(V,E)$ in which
every edge contains precisely $k$ distinct vertices of $V$. Quasi-random
graphs were first explicitly studied by Thomason \cite{T1,T2} and then followed by Chung, Graham, and
Wilson \cite{CGW}. Quasi-random properties were also studied in various other contexts
such as set systems \cite{CG}, tournaments \cite{CG1}, and
hypergraphs \cite{CG2}. There are also some very recent results on
quasi-random groups \cite{Go} and generalized quasi-random graphs
\cite{LS}. We briefly mention that the study of quasi-random
structures lies at the core of the recent proofs of Szemer\'edi's
Theorem \cite{Sztheo} that were recently obtained independently by
Gowers \cite{Go1,Go2} and by Nagle, R\"odl, Schacht and Skokan \cite{NRS,RS} and then
also by Tao \cite{T} and Ishigami \cite{Ish}. For more mathematical background on
quasi-randomness the reader is referred to the recent papers of
Gowers \cite{Go,Go1,Go2}. Quasi-random graphs are also related to theoretical computer-science via
the algorithmic version of the regularity lemma and the notion of expanders. For more details on quasi-random
graphs we refer the reader to the survey of Krivelevich and
Sudakov \cite{KS}.

We start with discussing quasi-random graphs.
One of the most natural questions that arise when studying quasi-random objects, is which properties
``force'' an object to behave like a truly-random one. The cornerstone result of this
type is the theorem of quasi-random graphs due to Chung, Graham and Wilson \cite{CGW}.
We start with some notation. For a subset of vertices $U$ in a graph $G$ we denote by $e(U)$ the number
of edges spanned by $U$ in $G$. For a pair of sets $U,U'$ we denote by $e(U,U')$ the number of edges with
one vertex in $U$ and the other in $U'$. Note that in a random graph $G(n,p)$ we expect every $U \subseteq V$
to satisfy $\frac12p|U|^2 -o(n^2) \leq e(U) \leq \frac12p|U|^2 + o(n^2)$, where here and throughout the paper an $o(1)$ term denotes (as usual)
any quantity that approaches $0$ as $n$ goes to infinity. To further simplify the notation, we will write $a=x\pm y$ to denote
the fact that $x-y \leq a \leq x+y$. So the above fact can be written as $e(U) = \frac12p|U|^2 \pm o(n^2)$.
The following is (part of) the main result of \cite{CGW}:

\begin{theo}[Chung-Graham-Wilson \cite{CGW}]\label{theoCGW}
Fix any $0 < p < 1$. For any $n$-vertex graph $G$ the following
properties are equivalent:
\begin{itemize}
\item ${\cal P}_1$: For any subset of vertices $U \subseteq
V(G)$ we have $~e(U)=\frac12p|U|^2 \pm o(n^2)$.
\item ${\cal P}_2(\alpha)$: For any subset of vertices $U \subseteq
V(G)$ of size $\alpha n$ we have $~e(U)=\frac12p|U|^2
\pm o(n^2)$.
\item ${\cal P}_3$: $e(G) = \frac{1}{2}p
n^2 \pm o(n^2)$ and $G$ has $\frac{1}{8}p^4n^4 \pm o(n^4)$ cycles of length $4$.
\end{itemize}
\end{theo}

As we have mentioned before, we use the $o(1)$ term to denote an arbitrary function tending to $0$ with
$n$. Hence, the meaning of the fact that, for example, ${\cal P}_2(1/2)$ {\em
implies} ${\cal P}_1$ is that for any $f(n)=o(1)$ there is a $g(n)=o(1)$
such that if $G$ has the property that
all $U \subseteq V(G)$ of size $n/2$ satisfy
$e(U)=\frac{1}{2}p|U|^2 \pm  g(n)n^2$, then
$e(U)=\frac{1}{2}p|U|^2 \pm f(n)n^2$ for all $U \subseteq V(G)$.
Equivalently, this means for any $\delta>0$ there is an
$\epsilon=\epsilon(\delta)$ and $n_0=n_0(\delta)$ such that if $G$ is
a graph on $n \geq n_0$ vertices and it has the property that
all $U \subseteq V(G)$ of size $n/2$ satisfy
$e(U)=\frac{1}{2}p|U|^2 \pm \epsilon n^2$, then
$e(U)=\frac{1}{2}p|U|^2 \pm \delta n^2$ for all $U \subseteq V(G)$.
This will also be the meaning of other implications between
other graph properties later on in the paper.


Note, that each of the items in Theorem \ref{theoCGW} is a property
we would expect $G(n,p)$ to satisfy with high probability. We will
thus say that $G$ is {\em $p$-quasi-random} if it satisfies property
${\cal P}_1$, that is if for some small $\delta$ all $U \subseteq
V(G)$ satisfy $e(U)=\frac{1}{2}p|U|^2 \pm \delta n^2$. If one wishes
to be more formal then one can in fact say that such a graph is
$(p,\delta)$-quasi-random. We will sometimes omit the $p$ and just
say that a graph is quasi-random. In the rest of the paper the
meaning of a statement ``If $G$ satisfies ${\cal P}_2$ then $G$ is
quasi-random'' is that ${\cal P}_2$ implies ${\cal P}_1$ in the
sense of Theorem \ref{theoCGW} discussed in the previous paragraph.
We will also say that a graph property ${\cal P}$ is {\em
quasi-random} if any graph that satisfies ${\cal P}$ must be
quasi-random. So the meaning of the statement ``${\cal P}_2$ is
quasi-random'' is that ${\cal P}_2$ implies ${\cal P}_1$. Therefore,
all the properties in Theorem \ref{theoCGW} are quasi-random.

As in Theorem \ref{theoCGW}, for the rest of the paper we fix $p$ to be the density of the ``supposed'' quasi-random
graph. Therefore all properties we will define from now on will refer to what one can expect to find in $G(n,p)$.

Our main focus in this paper is on the quasi-randomness of graph (and hypergraph) properties that involve the number of edges
in certain cuts in a graph (hypergraph). These properties were first studied by Chung and Graham \cite{CG,CG3}. We start by defining
the types of properties studied in \cite{CG,CG3}.

\begin{definition}\label{2cutgraph}
Fix an $0 < a < 1$. A graph satisfies property ${\cal P}_{a}$ if for any
$U \subseteq V(G)$ of size $|U|=an$ we have $e(U, V \setminus
U)=pa(1-a)n^2 +o(n^2)$.
\end{definition}

The main result of \cite{CG3,CG} was a precise characterization of the cut properties ${\cal P}_{a}$ which
are quasi-random. The somewhat surprising characterization is the following:

\begin{theo}[Chung-Graham \cite{CG,CG3}]\label{CGtheo} $P_{a}$ is quasi-random if and only if $a \neq 1/2$.
\end{theo}

To see that ${\cal P}_{1/2}$ is not quasi-random, Chung and Graham \cite{CG} observed that the graph obtained by taking
a random graph $G(n/2,2p)$ on $n/2$ of the vertices, an independent set on the other $n/2$ vertices and then connecting these
two graphs with a random bipartite graph with edge probability $p$ gives a non-quasi-random graph that satisfies ${\cal P}_{1/2}$.
For later reference, we call this graph $C_2(n,p)$. Chung and Graham \cite{CG} then
gave two proofs of the fact that when $a \neq 1/2$ property $P_{a}$ is quasi-random. One proof was based on a counting argument and another based on the rank of certain intersection matrices. Another proof of Theorem \ref{CGtheo}, using the machinery of graph limits, was given recently by Janson \cite{J}.

One of the open problems raised by Chung and Graham in their paper on quasi-random hypergraphs, was if one can obtain an analog of Theorem \ref{CGtheo} for $k$-uniform hypergraphs (see Section 9 in \cite{CG}).
Our first result in this paper answers this question positively by obtaining a precise characterization of the hypergraph cut properties that are quasi-random. This result
is discussed in Subsection \ref{subseccharac}. Our second result in this paper will show that one can ``describe'' the graphs (and hypergraphs) that satisfy the cut properties ${\cal P}_{a}$ which are not quasi-random. In particular, it will turn out that the example of Chung-Graham \cite{CG} (the graph $C_2(n,p)$ described above) showing that $P_{1/2}$ is not quasi-random is (essentially) the only graph that satisfies ${\cal P}_{a}$ and is not quasi-random. This result is discussed in Subsection \ref{subsecstruct}.

\subsection{A characterization of the quasi-random hypergraph cut properties}\label{subseccharac}

We now turn to discuss our first result which extends the result of Chung and Graham, stated in Theorem \ref{CGtheo},
from graphs to hypergraphs.
We will use the same notation we have used in the previous subsection for graphs. Let us first define a property of $k$-uniform
hypergraphs which is analogous to property ${\cal P}_1$ in Theorem \ref{theoCGW}.

\begin{definition}[${\cal D}_1$]
A $k$-uniform hypergraph $H=(V,E)$ satisfies property ${\cal D}_1$ if
$$
e(U)=\frac{p}{k!}|U|^k \pm o(n^k)
$$
for any $U \subseteq V$.
\end{definition}

So property ${\cal D}_1$ is perhaps the most intuitive notion of what it means for a hypergraph to be quasi-random. However, it turns out that
in many cases ${\cal D}_1$ is not the ``right'' generalization of ${\cal P}_1$. For example, while ${\cal P}_1$ implies that a graph has
the correct number of copies of any fixed graph (see \cite{CGW}), property ${\cal D}_1$ does not imply this fact. That is, for any $k>2$, there are $k$-uniform hypergraphs satisfying ${\cal D}_1$ that have {\em no} copy of (say) a clique of size $k+1$. This is the reason ${\cal D}_1$ is sometimes referred to as {\em weak} hypergraph quasi-randomness. There are more demanding notions of hypergraph quasi-randomness which do guarantee that
a quasi-random hypergraph would satisfy all the properties analogous to those guaranteed by ${\cal P}_1$.
This stronger notion of quasi-randomness was first defined by Frankl and R\"odl \cite{FR} and was recently studied due to its relation
to Szemer\'edi's Theorem \cite{Sztheo} in \cite{Go1,Go2,Ish,NRS,RS,T}. For an excellent discussion on the distinction between these notions of quasi-randomness, the reader is referred to \cite{Go1}. See
also \cite{CHPS,KNRS,NPRS} for other investigations of the notion of weak quasi-randomness in hypergraphs. We finally note that the reason
why the notion of weak quasi-randomness is interesting on its own, is that many properties are equivalent to weak quasi-randomness but not to the stronger notions. In particular, the properties we will study here will only be equal to weak quasi-randomness. Furthermore, while any hypergraph has a partition which is weakly quasi-random (in the sense of Szemer\'edi's regularity lemma \cite{Sz}) it is not true that any hypergraph has a strong quasi-random partition. Another important feature of weak hypergraph quasi-randomness is that it can be used as a tool to prove results about quasi-random graphs, as we shall also demonstrate in this paper. So henceforth, we will say that a hypergraph is $p$-quasi-random (or just quasi-random) if it satisfies ${\cal D}_1$ defined above.

We now define the hypergraph cut properties we will study in this paper.

\begin{definition}[${\cal P}_{\alpha}$]\label{cutprop}
Let $\alpha=(\alpha_1,\ldots,\alpha_r)$ be a vector of positive reals satisfying $\sum_i \alpha_i=1$.
For some $k \leq r$, we say that a $k$-uniform hypergraph on $n$ vertices satisfies property ${\cal P}_{\alpha}$ if for any partition of its vertices
into $r$-sets $V_1,\ldots,V_r$, where $|V_i|=\alpha_i n$, we have
$$
e(V_1,\ldots,V_r)=(p+o(1))n^k\sum_{S \subseteq [r], |S|=k} \,\prod_{i \in S}\alpha_i\;.
$$
\end{definition}

Here $e(V_1,\ldots,V_r)$ denotes the number of edges that cross the cut $(V_1,\ldots,V_r)$ (that is, the number of edges that intersect each $V_i$ in at most one
point). Note that this is a generalization of Definition \ref{2cutgraph} in that
we consider arbitrary $k$-uniform hypergraphs and in that we allow an arbitrary number of partition classes $r$ in the cut.
Chung and Graham \cite{CG} asked whether one can extend Theorem \ref{CGtheo} by finding characterization of the properties
${\cal P}_{\alpha}$ which are quasi-random, that is, equivalent to ${\cal D}_1$. Our first result in this paper answers this question by proving the following.

\begin{theo}\label{theomain}
Property ${\cal P}_{\alpha}$ is equivalent to ${\cal D}_1$ if and only if $\alpha \neq (1/r,\ldots,1/r)$.
\end{theo}

Given the graph $C_2(n,p)$ (defined after the statement of Theorem \ref{CGtheo}) which shows that ${\cal P}_{1/2}$ is not quasi-random, it seems natural
to try and show that when $\alpha=(1/r,\ldots,1/r)$ property ${\cal P}_{\alpha}$ is not quasi-random for $k$-uniform hypergraph by defining an appropriate
$k$-partite $k$-uniform hypergraph. This approach does not seem to work. Instead, we define the following $k$-uniform hypergraph.

\begin{definition}[$C_k(n,p)$]\label{defCnp} Let $C_k(n,p)$ be the $n$-vertex hypergraph constructed randomly as follows.
We partition the vertex set into two sets $A,B$ of size $n/2$ each. Each set of $k$ vertices $\{v_{i_1},\ldots,v_{i_k}\}$
is put in $C_k(n,p)$ with probability $2pj/k$ where $j=|\{v_{i_1},\ldots,v_{i_k}\}\cap A|$.
\end{definition}

Observe that when $k=2$ the graph $C_k(n,p)$ defined above is (indeed) equivalent to the (randomly constructed) graph $C_2(n,p)$ we described earlier.
As we will show later, this random hypergraph satisfies ${\cal P}_{\alpha}$ (for $\alpha=(1/r,\ldots,1/r)$) with high probability
but is not quasi-random, that is, does not satisfy ${\cal D}_1$ defined above. This will establish that ${\cal P}_{\alpha}$ is not quasi-random. Our second result in this paper,
discussed in the next subsection, shows that the hypergraphs $C_k(n,p)$ are essentially
the only non quasi-random hypergraphs satisfying ${\cal P}_{\alpha}$.

\subsection{The structure of graphs satisfying a non quasi-random cut property}\label{subsecstruct}

The fact that a graph property ${\cal P}$ is quasi-random means that knowing that a graph satisfies ${\cal P}$ tells
us a lot about the structure of the graph. It is natural to ask if knowing that a graph satisfies a {\em non} quasi-random
graph property, still tells us something about its structure \footnote{Of course, this question is not well defined but the spirit of
it should be clear to a reader who is familiar with the notion of quasi-random graphs.}. For example, while we learn from Theorem \ref{theoCGW} that a graph with the correct number
of edges and the correct number of copies of $C_4$ (the $4$-cycle) must be quasi-random, this is no longer the case if one considers the number of triangles
rather then the number of copies of $C_4$. Furthermore, it does not seem like one can ``describe'' the graphs that have the correct number of edges and the correct
number of $K_3$. Our second main result is that when considering the non quasi-random cut properties then one can obtain such a description.

Let's consider first the case of graphs. In this case the non-quasi-random cut property is ${\cal P}_{1/2}$ which corresponds to counting the number of edges in balanced $(n/2,n/2)$-cuts.
To describe our structure result about the graphs satisfying ${\cal P}_{1/2}$ it will be more convenient to consider the following non-discrete version of ${\cal P}_{1/2}$ which
we denote ${\cal P}^{*}_{1/2}$; in this problem we are asked to assign arbitrary real weights to the edges of the complete graph on $n$ vertices in a way that for any partition
of its vertices into two sets of equal size $n/2$, the total weight of edges crossing the cut is $p(n/2)^2$. Note that since ${\cal P}^{*}_{1/2}$ allows for non-integer weights, we
require the total weight crossing the cuts to be {\em exactly} $p(n/2)^2$, while in ${\cal P}_{1/2}$ the requirement is only up to an error of $o(n^2)$.

Considering the fractional property ${\cal P}^{*}_{1/2}$ we now ask which weight assignments satisfy ${\cal P}^{*}_{1/2}$? Observe that this problem can be stated as trying to solve a set of linear
equations, where for every $i<j$ we have an unknown $x_{i,j}$ and where for every partition of the $n$ vertices into two sets of equal size $n/2$, we have a linear equation $\ell_{A,B}$
which checks whether $\sum_{i \in A,j\in B}x_{i,j}=p(n/2)^2$. So this set has ${n \choose 2}$ unknowns and ${n-1 \choose n/2-1}$ equations.
One solution to this set of equations is the one corresponding to the random graph
$G(n,p)$ in which all $x_{i,j}=p$. Another solution corresponds to the graph $C_2(n,p)$ from Definition \ref{defCnp}.
In this case, we obtain a solution by partitioning the vertices into two sets $A$ and $B$ of size $n/2$ each, and setting $x_{i,j}=2p$ if $i,j \in A$, setting $x_{i,j}=0$ if $i,j \in B$ and setting $x_{i,j}=p$ otherwise.
Note that we thus obtain ${n-1 \choose n/2-1}$ solutions which correspond to the possible ways of picking the sets $A,B$. However, observe that all these solutions are isomorphic to $C_2(n,p)$, if we consider them as weighted complete graphs.

So we can restate our question and ask if there are any other solutions to ${\cal P}^{*}_{1/2}$ besides the above $1+{n-1 \choose n/2-1}$ solutions?
Since we are trying to solve a set of linear equations, then one can trivially obtain other solutions by taking affine combinations of the above solutions. That is, if one considers each of the above
solutions as an ${n \choose 2}$ dimensional vector, then any affine combination of these vectors is also a solution. Our second result in this paper states that these are the only solutions to ${\cal P}^{*}_{1/2}$.

\begin{theo}\label{theo2-informal}
The only solutions to ${\cal P}^{*}_{1/2}$ are the affine combinations of $G(n,p)$ and $C_2(n,p)$.
\end{theo}

So the above theorem can be restated as saying that the only graphs satisfying ${\cal P}_{1/2}$ are those that can be obtained in a trivial way from the random graph $G(n,p)$ and the counter example
of Chung-Graham showing that ${\cal P}_{1/2}$ is not quasi-random. As we show in Section \ref{secstruct}, given Theorem \ref{theo2-informal} one can easily show that any graph satisfying ${\cal P}_{1/2}$ can be approximated by an affine combination of $G(n,p)$ and $C_2(n,p)$, thus
supplying a structural characterization of the graphs satisfying ${\cal P}_{1/2}$. See Theorem \ref{cutnorm} in Section \ref{secstruct}.

When considering the hypergraph cut properties ${\cal P}_{\alpha}$ of Definition \ref{cutprop}, we can of course define
${\cal P}^{*}_{\alpha}$ to be their non-discrete analog. That is, we now try to assign weights to the edges of the
complete $k$-uniform hypergraph. In Section \ref{secstruct} we also prove the following theorem which extends Theorem \ref{theo2-informal} to hypergraphs. In the following statement $G_k(n,p)$ denotes the random $k$-uniform hypergraph on $n$ vertices.

\begin{theo}\label{theo3-informal}
Set $\alpha=(1/r,\ldots,1/r)$. The only solutions to ${\cal P}^{*}_{\alpha}$ are the affine combinations of $G_k(n,p)$ and $C_k(n,p)$.
\end{theo}


\subsection{Organization}

The rest of the paper is organized as follows. In Section \ref{Secunbalanced} we prove the first implication of Theorem
\ref{theomain} showing that unbalanced cuts are quasi-random. The proof has two main steps. In the first step we reprove
the result of Chung and Graham \cite{CG} on cuts in graphs using a simple argument, which uses a method that was recently introduced
by the authors in \cite{S,SY,Y} for tackling problems on quasi-random graphs and hypergraphs. This method uses probabilistic and algebraic arguments to analyze the edge distribution in graphs and hypergraphs. See also \cite{CHPS} where this method was used to study weak quasi-random
hypergraphs. We then prove Theorem \ref{theomain} by reducing it to the argument we use in order to reprove Theorem \ref{CGtheo}.
The other side of Theorem \ref{theomain} is proved in Section \ref{seccounter} where we prove that the non quasi-random hypergraphs $C_k(n,p)$ defined above satisfy
${\cal P}_{\alpha}$ when $\alpha$ is a balanced cut.
Theorems \ref{theo2-informal} and \ref{theo3-informal} are proved in Section \ref{secstruct}.
Both proofs rely on the computation of the rank of certain intersection matrices which we carry out in Section \ref{Secmatrix}. While the proof of Theorem \ref{theo2-informal} is technically simple,
the proof of Theorem \ref{theo3-informal} is much more involved and applies certain results from the theory of Association Schemes.
We believe the ideas here may be applicable to the study of other properties of quasi-random graphs and hypergraphs.
Finally, in Section \ref{secopen} we discuss another generalization of Theorem \ref{CGtheo} and raise a related open problem.

\section{Unbalanced Cut Properties Are Quasi-Random}\label{Secunbalanced}

We begin this section with the first implication of Theorem \ref{theomain}. Actually, we prove that property ${\cal D}_1$ implies property
${\cal P}_{\alpha}$ for any $\alpha$ (whether balanced or not).

\begin{lemma}\label{easy} If a $k$-uniform hypergraph $H$ satisfies ${\cal D}_1$ then for any $\alpha=(\alpha_1,\ldots,\alpha_r)$ it also satisfies
$P_{\alpha}$.
\end{lemma}

\paragraph{Proof:} We start with graphs. Let $V_1,\ldots,V_r$ be a partition of $V(H)$ into $r$ sets satisfying $|V_i|=\alpha_in$. Since $H$ satisfies ${\cal D}_1$ we have $e(U)=\frac12p|U|^2 \pm o(n^2)$ for all $U \subseteq V(H)$. Therefore, for any pair of disjoint sets $A$ and $B$ of sizes $\alpha n$ and $\beta n$ we have
\begin{eqnarray*}
e(A,B)&=&e(A \cup B)-e(A)-e(B)\\
&=&\frac12p(\alpha+\beta)^2n^2 \pm o(n^2)-\frac12p\alpha^2n^2 \pm o(n^2)-\frac12p\beta^2n^2 \pm o(n^2)\\
&=& \alpha\beta pn^2 \pm o(n^2)\;.
\end{eqnarray*}
Summing over all pairs $V_i,V_j$ we get $e(V_1,\ldots,V_r)=pn^2 \sum_{i < j}\alpha_i \alpha_j \pm o(n^2)$, as needed.

As to $k$-uniform hypergraphs, if we take $k$ vertex sets of sizes $\alpha_1n,\ldots,\alpha_kn$ then since every set of vertices $U$ spans $\frac{p}{k!}|U|^k \pm o(n^k)$, then by Inclusion-Exclusion we have
\begin{eqnarray*}
e(V_1,\ldots,V_k)&=& \sum^1_{t=k}(-1)^{k-t}\left(\sum_{S \subseteq [k]:|S|=t} e(\bigcup_{i \in S}V_i)\right)\\
&=& \sum^1_{t=k}(-1)^{k-t}\left(\sum_{S \subseteq [k]:|S|=t}  \frac{p}{k!}\left(\sum_{i \in S}\alpha_i\right)^kn^k \pm o(n^k) \right)\\
&=& \frac{p}{k!}n^k \sum^1_{t=k}(-1)^{k-t}\left(\sum_{S \subseteq [k]:|S|=t} \left(\sum_{i \in S}\alpha_i\right)^k  \right)\pm o(n^k)\\
&=& pn^k \prod^k_{i=1}\alpha_i  \pm o(n^k)\;.
\end{eqnarray*}
The last equality follows from the observation that when expanding the expression
$$
\sum^1_{t=k}(-1)^{k-t}\left(\sum_{S \subseteq [k]:|S|=t} \left(\sum_{i \in S}\alpha_i\right)^k  \right)
$$
we get a sum of monomials of the form $\prod_{i \in S}\alpha^{x_i}_i$ with $\sum_{i \in S}x_i=k$.
The coefficient of such a monomial is
$$
\frac{k!}{x_1! \cdots x_k!}\sum^{|S|}_{t=k}(-1)^{k-t}{k-|S| \choose t-|S|}
$$
which is $0$ when $|S| < k$ and $k!$ when $|S|=k$ (in this case $x_1=\ldots=x_k=1$). We now complete the proof as in the case of graphs by summing over all collections of $k$ subsets $V_{i_1},\ldots,V_{i_k}$.
$\qed$

\bigskip

We now turn to the proof of the second (and main) implication of Theorem \ref{theomain}.
As we have mentioned earlier, we will first give a simple and short proof
of the main result of Chung and Graham \cite{CG}, stated in Theorem \ref{CGtheo}, which deals with the special case of
graphs. We will then prove the general case by reducing it to the case $k=2$. Let us restate the result of \cite{CG}.

\begin{theo}[\cite{CG}]\label{graphk2}
Let $\alpha=(a,1-a)$, where $0 < a < 1$. If $a \neq 1/2$ then ${\cal P}_{\alpha}$ is quasi-random.
\end{theo}

\paragraph{Proof:} We will prove that if $a \neq 1/2$ then ${\cal P}_{\alpha}$ implies ${\cal P}_2(a)$
and is thus $p$-quasi-random by Theorem \ref{theoCGW}. Clearly we may assume that $a < 1/2$.
So fix any set $A$ of $an$ vertices and let $B=V-A$. Let $x_0,x_1,x_2$ satisfy
$|e(A)|=x_{0}\frac12a^2 n^2$, $|e(A,B)|=x_{1}a(1-a)n^2$ and $|e(B)|=x_{2}\frac12(1-a)^2n^2$.
We need to show that $x_0=p \pm o(1)$.

Let $0 \leq c \leq a$ and consider a (new) vertex partition $V_1,V_2$ of $G$ into sizes $an,(1-a)n$
that is constructed as follows: We randomly select $cn$ vertices of $A$ and place them in $V_1$, and randomly
select $(a-c)n$ vertices of $B$ and also place them in $V_1$.
The probability that an edge belonging to $A$ connects $V_1$ to $V_2$ is $2{an-2 \choose cn-1}/{an \choose cn}=\frac{2c(a-c)}{a^2} \pm o(1)$. Similarly, the probability that
an edge with one endpoint in $A$ and another in $B$ connects $V_1$ to $V_2$ is $\frac{c(1-2a+c)+(a-c)^2}{a(1-a)}$,
and the probability that an edge belonging to $B$ connects $V_1$ to $V_2$ is $2{(1-a)n-2 \choose (a-c)n-1}/{(1-a)n \choose (a-c)n}=\frac{2(a-c)(1-2a+c)}{(1-a)^2}\pm o(1)$.
Therefore, the expected number of edges connecting $V_1$ to $V_2$ is
$$
c(a-c)x_{0} n^2 + (c(1-2a+c) + (a-c)^2)x_{1}n^2 + (a-c)(1-2a+c)x_{2}n^2 \pm o(n^2)\;.
$$
But since we assume that $e(V_1,V_2)=a(1-a)pn^2 \pm o(n^2)$ for {\em every} $(a,1-a)$-cut in $G$, this expectation must
equal $a(1-a)pn^2 \pm o(n^2)$. Dividing by $n^2$ we get that for every $0 \leq c \leq a$
\begin{equation}\label{k2}
c(a-c)x_{0} + (c(1-2a+c) + (a-c)^2)x_{1} + (a-c)(1-2a+c)x_{2} = a(1-a)p \pm o(1)\;.
\end{equation}

Assume for a moment that (\ref{k2}) holds without the $o(1)$ term. Using the values
$c=0$, $c=a$ and $c=a/2$ we get three linear equations $Ax=a(1-a)p$ where
$A$ is the matrix
$$
\left(
\begin{array}{ccc}
0 &  a^2 & a(1-2a) \\
0 & a(1-a) & 0\\
a^2/4 & a(1-a)/2 & a(1-3a/2)/2
\end{array}
\right)
$$
Since $Det(A)=\frac{a^2}{4}a(1-a)a(1-2a)$ and we assume that $a \not \in \{0,\frac12,1\}$ we have $Det(A) \neq 0$ and so
$Ax=a(1-a)p$ has a unique solution. As $x_0=x_1=x_{2}=p$ is a valid solution of $Ax=a(1-a)p$, this {\em is} the (unique) solution. Since a solution of a system of linear equation $Ax=b$ is continuous with respect to $b$, we get
that when considering (\ref{k2}) with the $o(1)$ term, the solution still satisfies $x_0,x_1,x_2 = p \pm o(1)$, thus completing the proof. $\qed$

\bigskip

For the second part of the proof of Theorem \ref{theomain}, we will need to consider a relaxed version of property ${\cal D}_1$, analogous to property ${\cal P}_2$ in Theorem \ref{theoCGW}. We will need the following lemma.

\begin{lemma}\label{P2hyper} Fix $0 < \eta < 1$ and let ${\cal D}_{\eta}$ be the hypergraph property of satisfying
$e(U)=\frac{p}{k!}|U|^k \pm o(n^k)$ for all $U \subseteq V$ of size $\eta n$. Then properties ${\cal D}_{\eta}$ and
${\cal D}_1$ are equivalent.
\end{lemma}

In \cite{Y} the second named author proves that if every $\eta n$ vertices in a graph
contain the ``correct'' number of copies of $K_k$ one expects to find in $G(n,p)$ then {\em every} set contains
the correct number of copies of $K_k$. Precisely the same proof gives the above lemma. Hence we refrain from reproducing the identical proof.

\begin{theo}\label{hypergraphkk}
Let $\alpha=(\alpha_1,\ldots,\alpha_k)$ be a vector of positive reals satisfying $\sum_i \alpha_i=1$.
If $\alpha \neq (1/k,\ldots,1/k)$ then ${\cal P}_{\alpha}$ is quasi-random.
\end{theo}

\paragraph{Proof:} The case $k=2$ was handled in Theorem \ref{graphk2}, so we assume that $k >2$.
Suppose without loss of generality that $\alpha_{k-1} < \alpha_k$.
We will prove that every set $A$ of $(1-\alpha_k)n$ vertices has
$\frac{p}{k!}|A|^k \pm o(n^k)$ edges, and hence the result follows by Lemma \ref{P2hyper}.
So fix a set $A$ of size $(1-\alpha_k)n$ and let $B=V-A$.
For $0 \leq i \leq 2$ let $x_i$ denote the density of edges with $i$ vertices in $B$ and $k-i$ vertices in $A$.
We need to show that $x_0=p \pm o(1)$, but as in Theorem \ref{k2} it will be easier
to show that all three densities are $p \pm o(1)$.

Take any $0 \leq d \leq \alpha_{k-1}$ and consider a vertex partition $V_1,V_2,\ldots,V_k$ of $H$ into $k$ sets of sizes $\alpha_in$ for $1 \leq i \leq k$,
that is constructed as follows. We randomly select
$\alpha_1n$ vertices of $A$ and put them in $V_1$. We then select $\alpha_2n$ of the remaining vertices of $A$ and put then in $V_2$.
We continue in the same manner by selecting $\alpha_in$ vertices of $A$ and put them in $V_i$ for $i=1,\ldots,k-2$.
We then randomly select $dn$ vertices of the remaining vertices of $A$ to place in $V_{k-1}$,
and randomly select $(\alpha_{k-1}-d)n$ vertices of $B$ to place in $V_{k-1}$.
The remaining vertices of $A$ (there are $(\alpha_{k-1}-d)n$ such vertices) are placed in $V_k$,
and so are the remaining vertices of $B$ (there are $(\alpha_k-\alpha_{k-1}+d)n$ such vertices).
The probability that an edge in $A$ has one point in each of the sets $V_1,\ldots,V_k$ is
$$
p_1=k!\left(\prod_{i=1}^{k-2}\frac{\alpha_i}{|A|/n}\right) \frac{d(\alpha_{k-1}-d)}{(|A|/n)^2} + o(1)\;.
$$
Similarly, the probability that
an edge with one point in $B$ and $k-1$ points in $A$ has one vertex in each of the sets $V_1,\ldots,V_k$ is
$$
p_2=(k-1)! \left(\prod_{i=1}^{k-2}\frac{\alpha_i}{|A|/n}\right) \frac{d(\alpha_k-\alpha_{k-1}+d)+(\alpha_{k-1}-d)^2}{|A||B|/n^2}+o(1)\;,
$$
and the probability that an edge with two points in $B$ and $k-2$ in $A$ has one point in each of the sets $V_1,\ldots,V_k$ is
$$
p_3=2(k-2)!\left(\prod_{i=1}^{k-2}\frac{\alpha_i}{|A|/n}\right) \frac{(\alpha_{k-1}-d)(\alpha_k-\alpha_{k-1}+d)}{(|B|/n)^2}+o(1)\;.
$$

As in the proof of Theorem \ref{k2} we get that the expected number of edges crossing $(V_1,\ldots,V_k)$ is
$$
p_1{{|A|} \choose k}x_0 + p_2{{|A|} \choose {k-1}}|B|x_1 + p_3{{|A|} \choose {k-2}}{{|B|} \choose 2}x_2
$$
and by the assumed properties of $H$ this expectation should be equal to $n^k\left(\prod_{i=1}^{k}\alpha_i\right)p \pm o(n^k)$.
Dividing by $n^k$ we get for every $0 \leq d \leq \alpha_{k-1}$ the equation
\begin{eqnarray*}
&&\left(\prod_{i=1}^{k-2}\alpha_i\right)d(\alpha_{k-1}-d)x_0+\\
&&\left(\prod_{i=1}^{k-2}\alpha_i\right)(d(\alpha_k-\alpha_{k-1}+d)+(\alpha_{k-1}-d)^2)x_1+\\
&&\left(\prod_{i=1}^{k-2}\alpha_i\right)(\alpha_{k-1}-d)(\alpha_k-\alpha_{k-1}+d)x_2
= \left(\prod_{i=1}^{k}\alpha_i\right)p \pm o(1)\;.
\end{eqnarray*}

Dividing the above equation by $\left(\prod_{i=1}^{k-2}\alpha_i\right)$ we get
$$
d(\alpha_{k-1}-d)x_0+(d(\alpha_k-\alpha_{k-1}+d)+(\alpha_{k-1}-d)^2)x_1+
$$
$$
(\alpha_{k-1}-d)(\alpha_k-\alpha_{k-1}+d)x_2 = \alpha_{k-1}\alpha_kp \pm o(1)\;.
$$

Let us now make a syntactic change in the variables. Define:
$r=\alpha_{k-1}+\alpha_k$, $c=d/r$ and $a=\alpha_{k-1}/r$.
We can rewrite the last equality as:
$$
r^2c(a-c)x_0 + r^2(c(1-2a+c)+(a-c)^2)x_1 + r^2(a-c)(1-2a+c)x_2 = \alpha_{k-1}\alpha_kp \pm o(1)\;.
$$
Dividing everything by $r^2$ we obtain for every $0 \leq c \leq  a < 1/2$ the linear equation
$$
c(a-c)x_0 + (c(1-2a+c)+(a-c)^2)x_1 + (a-c)(1-2a+c)x_2 = a(1-a)p \pm o(1)\;.
$$
This is exactly the same equations we got in the case $k=2$ in (\ref{k2}). By the argument given in the proof of Theorem \ref{k2},
we get that $x_0=p \pm o(1)$ which is what we wanted to prove. $\qed$

\bigskip

We now turn to extend the above result to cut properties with an arbitrary number of classes.

\begin{theo}\label{graphk}
Let $\alpha=(\alpha_1,\ldots,\alpha_r)$ be a vector of positive reals satisfying $\sum_i \alpha_i=1$.
If $\alpha \neq (1/r,\ldots,1/r)$ then ${\cal P}_{\alpha}$ is quasi-random.
\end{theo}

\paragraph{Proof:} We prove the theorem for graphs. The proof for $k$-uniform hypergraphs is identical.
Suppose without loss of generality that
$\alpha_{r-1} \neq \alpha_r$ and consider any set $A$ of $(\alpha_{r-1}+\alpha_r)n$ vertices in $G$. Fix any partition of $V \setminus A$
into $r-2$ sets $V_1,\ldots,V_{r-2}$ of sizes $|V_i|=\alpha_i n$. As $G$ satisfies ${\cal P}_{\alpha}$ we get that for any partition of $A$ into two sets $A_1,A_2$ of sizes $\alpha_{r-1}n$ and $\alpha_rn$ we have
\begin{eqnarray*}
(p+o(1))n^2\prod_{1 \leq i < j \leq r}\alpha_i\alpha_j &=&e(V_1,\ldots,V_{r-2},A_1,A_2)\\
&=&e(V_1,\ldots,V_{r-2})+e(\bigcup^{r-2}_{i=1}V_i,A)+e(A_1,A_2)\;.
\end{eqnarray*}
Therefore, for every partition of $A$ into sets $A_1$, $A_2$ of sizes $\alpha_{r-1}n$ and $\alpha_rn$ we have
\begin{eqnarray*}
e(A_1,A_2)&=&(p+o(1))n^2\prod_{1 \leq i < j \leq r}\alpha_{i,j}-e(V_1,\ldots,V_{r-2})-e(\bigcup^{r-2}_{i=1}V_i,A)\\
&=& (q+o(1))n^2\;,
\end{eqnarray*}
for some $0 \leq q \leq 1$ (which is independent of the partition $A_1,A_2$). Since we assume that $\alpha_r \neq \alpha_{r-1}$ we deduce using
Theorem \ref{graphk2} that the graph induced by $A$ is $q/(\alpha_{r-1}\alpha_{r})$-quasi-random. We claim that this means that every set of vertices in $G$ of size $(\alpha_{r-1}+\alpha_r)n$ spans (asymptomatically) the same number of edges. Indeed, take any pair of sets $A$ and $B$ of size $(\alpha_{r-1}+\alpha_r)n$ and let $C$ be a set of size $(\alpha_{r-1}+\alpha_r)n$ containing at least half of the vertices of each of the sets $A$ and $B$. Since $A$, $B$ and $C$ all span quasi-random graphs and $|A \cap C| \geq |A|/2,|C|/2$ and $|A \cap B| \geq |B|/2,|C|/2$ we have
$$
d(A)=d(A \cap C) \pm o(1)= d(C) \pm o(1)= d(B \cap C) \pm o(1)=d(B) \pm o(1)\;.
$$
Therefore, by property ${\cal P}_2$ of Theorem \ref{theoCGW} we have that $G$ is $p'$-quasi-random for some $p'$. Finally, by Lemma \ref{easy} we know that this means that $G$ satisfies ${\cal P}_{\alpha}$ with edge density $p'$. But since
we assume that $G$ satisfies ${\cal P}_{\alpha}$ with edge density $p$, we get that $p'=p$ and so $G$ is $p$-quasi-random, as needed. $\qed$

\section{Balanced Cut Properties are Not Quasi-Random}\label{seccounter}

We prove that for each $r \ge k \ge 2$, a hypergraph obtained by the random construction $C_k(n,p)$ satisfies ${\cal P}_{\alpha}$ for $\alpha=(1/r,\ldots,1/r)$,
with high probability. Recall that this is known for $r=k=2$ \cite{CG} and we now generalize this to all $k$ and $r$.
Our construction will assume that $p \leq \frac12$ (although it is not difficult to modify the construction to accommodate the case $p > 1/2$).

Recall that $C_k(n,p)$ is obtained via the following random construction.
We partition the vertex set into two parts $A,B$ with $|A|=|B|=n/2$.
For any $k$-subset of vertices, we randomly and independently select it to be an edge according to the following rule;
if the set has $j$ vertices in $A$ then it will be an edge with probability $2pj/k$. In particular $A$ induces a hypergraph with expected density $2p$
and $B$ induces an empty hypergraph. This clearly means that w.h.p. $C_k(n,p)$ is not quasi-random.

We now prove that $C_k(n,p)$ has the property that in every balanced $r$-cut, the expected density of the cut edges is $p$.
Since the number of ways to partition the vertex set into $r$ sets is bounded by $2^{rn}$, while the probability that the number of edges
in a given cut significantly deviates from its expectation is $2^{-\Theta(n^k)}$ (via a standard Chernoff bound) we get from the
union bound that with high probability the resulting hypergraph $C_k(n,p)$ satisfies ${\cal P}_{\alpha}$.

Let ${\bf z}=(z_1,\ldots,z_r)$ be a vector of positive reals with
\begin{equation}\label{L1}
\sum_{i=1}^r z_i=r/2\;.
\end{equation}
A balanced $r$-cut of $A \cup B$ with parts $V_1,\ldots,V_r$ is of {\em type ${\bf z}$} if $V_i$ contains precisely
$nz_i/r$ vertices of $A$ for every $1 \leq i \leq r$.

We fix ${\bf z}$ and show that cuts of type ${\bf z}$ have expected density $p$. Since this will hold for each fixed ${\bf z}$,
this will also hold for all balanced $r$-cuts, as required.

So consider some cut of type ${\bf z}$. We can write a closed formula for the expected density of edges in the cut.
But before that, let us see a simplified example for the case $k=2$ and $r=3$.
In this case, the expected density is
\begin{eqnarray*}
&&\frac{1}{3}(
2p(z_1z_2+z_1z_3+z_2z_3)+\\
&&p(z_1(1-z_2)+(1-z_1)z_2+z_1(1-z_3)+(1-z_1)z_3+z_2(1-z_3)+(1-z_2)z_3))\;.
\end{eqnarray*}
It is straightforward to verify that the above expression is identically $p$, as required.

Doing the same for general $k$ and $r$ requires, however, more care, and some notation.
Let ${X \choose y}$ denote the set of all $y$-element subsets of a set $X$.
For $K \in {{[r]} \choose k}$ and for $J \in {K \choose j}$ denote by $z_{K,J}$ the polynomial expression
$$
z_{K,J} = \left(\prod_{j \in J}z_j\right)\cdot \left(\prod_{j \in K-J}(1-z_j)\right)\;.
$$
Now, let
$$
s_{K,j}=\sum_{J \in {K \choose j}}z_{K,J}\;
$$
and let
$$
s_{r,k,j}=\sum_{K \in {{[r]} \choose k}}s_{K,j}\;.
$$
Hence, for example, if $r=4$, $k=3$ and $j=2$ we have
$$
s_{4,3,2}=
z_1z_2(1-z_3)+z_1(1-z_2)z_3+(1-z_1)z_2z_3+
z_1z_2(1-z_4)+z_1(1-z_2)z_4+(1-z_1)z_2z_4+
$$
$$
z_1z_3(1-z_4)+z_1(1-z_3)z_4+(1-z_1)z_3z_4+
z_2z_3(1-z_4)+z_2(1-z_3)z_4+(1-z_2)z_3z_4\;.
$$
In general, the density of the edges of cuts of type ${\bf z}$ is
\begin{equation}
\label{z}
\frac{1}{{r \choose k}}\left(\sum_{j=1}^r  \frac{2pj}{k} \cdot s_{r,k,j}\right)\;.
\end{equation}
It is therefore our goal to prove that
\begin{equation}\label{goal}
\sum_{j=1}^r  (2pj/k)s_{r,k,j} \equiv {r \choose k}p\;.
\end{equation}

When expanding (\ref{z}) we obtain a polynomial in $z_1,\ldots,z_k$. This polynomial is a sum of
monomials where each monomial is of the form $c_J\prod_{j \in J}z_j$ for some $J \in {{[r]} \choose j}$,
with $|J| \le k$, where $c_J$ is some constant. We will prove that $c_J=0$ for every $|J| > 1$.
Let us first examine the case where $J$ has $k$ elements.
First notice that $s_{r,k,k}$ contributes $2p$ to the constant $c_{J}$.
When expanding $s_{r,k,k-1}$ we notice that $\prod_{j \in J}z_j$ appears, with negative sign, $k$ times, and hence
$s_{r,k,k-1}$ contributes $-(2p(k-1)/k)\cdot k = -2p(k-1)$ to the constant $c_{J}$.
More generally, when expanding $s_{r,k,i}$ we notice that $\prod_{j \in J}z_j$ appears ${k \choose i}$
times and with sign $(-1)^{k-i}$. Hence we obtain
$$
c_{J} = \sum_{i=1}^k (-1)^{k-i}\frac{2pi}{k} {k \choose i} = 2p \sum_{i=0}^{k-1} (-1)^{k-1-i} {{k-1} \choose {i}} \equiv 0.\;
$$

Now let us examine a general $J$ with $|J| > 1$. Set $|J|=j$. We first notice that $s_{r,k,q}$ does not contribute anything to
$c_J$ whenever $q > j$. When expanding $s_{r,k,j}$ we notice that $\prod_{j \in J}z_j$ appears, with positive sign, precisely once.
More generally, when expanding $s_{r,k,i}$ for $1 \le i \le j$ we notice that $\prod_{j \in J}z_j$ appears ${j \choose i}$ times
and with sign $(-1)^{j-i}$. Hence we obtain
$$
c_{J} = \sum_{i=1}^j (-1)^{j-i}\frac{2pi}{k} {j \choose i} = \frac{2pj}{k} \sum_{i=0}^{j-1} (-1)^{j-1-i} {{j-1} \choose {i}} \equiv 0.\;
$$

It remains to consider the coefficients of the singletons $z_i$ for $i=1,\ldots,r$. Clearly, only the expansion of $s_{r,k,1}$ contains
the singleton $z_i$, and precisely ${{r-1} \choose{k-1}}$ times. Hence, it follows from (\ref{L1}) that
$$
\sum_{j=1}^r  (2pj/k)s_{r,k,j} = (2p/k)s_{r,k,1} = {{r-1} \choose{k-1}}\frac{2p}{k}\sum^k_{i=1}z_i \equiv {r \choose k}p\;,
$$
thus verifying (\ref{goal}) and completing the proof.

\section{The Rank of Certain Intersection Matrices}\label{Secmatrix}

Intersection matrices have been extensively studied for many years, see e.g. \cite{BF}. The key ingredient we need for the proofs of Theorems \ref{theo2-informal} and \ref{theo3-informal} regarding
the structure of graphs and hypergraphs which satisfy the balanced cut properties is Theorem \ref{RankA} below.
This theorem determines the exact rank of a certain intersection matrix. We start with formally defining the types of matrices we study.

\begin{definition}[$A_{t,k,v}$]\label{DefMatrix} For a positive vector $v=(v_1,\ldots,v_k)$ satisfying $\sum_iv_i=t$ let
$A=A_{t,k,v}$ be the following $0/1$ matrix. The columns of $A$ are indexed by the subsets of $\{1,\ldots,t\}$ of size $k$ and the rows
are indexed by the partitions of $\{1,\ldots,t\}$ into $k$ sets of sizes $v_1,\ldots,v_k$. With this indexing of the rows and columns, we set $A_{i,j}=1$ if and only if the $k$-set
$S$ corresponding to index $j$ has exactly one element in each of the sets of the partition $V_1,\ldots,V_k$ whose index is $i$.
\end{definition}

Since $A_{t,k,v}$ has ${t \choose k}$ columns we trivially have $rank(A_{t,k,v}) \leq {t \choose k}$. The following theorem gives a precise bound
for the rank of $A_{t,k,v}$. To avoid degenerate cases (where the number of rows of $A_{t,k,v}$ is smaller than the number of columns) we only
consider vectors $v$ where each coordinate is at least $k$.

\begin{theo}\label{RankA} For every $k$ and large enough $t \geq t_0(k)$, the following holds for every $v$:
\begin{equation}
rank(A_{t,k,v})=
\left\{%
\begin{array}{ll}
    {t \choose k}-t+1 & \hbox{$v=(t/k,\ldots,t/k)$} \\
    {t \choose k} & \hbox{$v\neq(t/k,\ldots,t/k)$} \\
\end{array}%
\right.
\end{equation}
\end{theo}

We note that Theorem \ref{RankA} can be used to give an alternative proof of Theorem \ref{theomain}. However, since the proof
of Theorem \ref{RankA} is much more complicated than the proof of Theorem \ref{theomain} we decided to give the more elementary proof described
in Section \ref{Secunbalanced}. Also, the case $k=2$ and $v \neq (t/2,t/2)$ was already considered by Chung and Graham \cite{CG},
who used the fact that in this case $A_{t,2,v}={t \choose 2}$ in order to give one of their proofs
of Theorem \ref{CGtheo}. As we show at the end of this section, one can actually prove this special case of Theorem \ref{RankA} by a direct reduction
to Gotllieb's Theorem \cite{G}. Actually, the same proof will work for all $k$ when $v \neq (t/k,\ldots,t/k)$, see Lemma \ref{fullrank}.

The hardest (and more interesting) part of proving Theorem \ref{RankA} is the case $v=(t/k,\ldots,t/k)$. It will actually be easier to obtain
the fact that in this case $rank(A_{t,k,v}) \leq {t \choose k}-t+1$ as part of our discussion in the next section, see Lemma \ref{l-rank3}.
We are thus left with proving that when $v=(t/k,\ldots,t/k)$ we have $rank(A_{t,k,v}) \geq {t \choose k}-t+1$. As we explain at the end of this section,
the special case $k=2$ is relatively easy to handle due to a certain degeneracy of this case.
Indeed the proof for $k \geq 3$ is much more complicated and is the main focus of this section. The proof will apply certain
results from the theory of association schemes discussed below. Hence, we now turn to prove the following lemma.

\begin{lemma}\label{hardpart} For every $k \geq 2$ and large enough $t \geq t_0(k)$, if
$v=(t/k,\ldots,t/k)$ then $rank(A_{t,k,v}) \geq {t \choose k}-t+1$.
\end{lemma}

Before getting to the details of the proof, we need to introduce some concepts from the theory of Association Schemes.
For integers $t \ge k \ge 2$ we define a set of $k+1$ symmetric binary matrices $J(t,k)=\{W_0,\ldots,W_k\}$
as follows. The rows and columns of each of the $W_i$ are indexed by ${{[t]} \choose k}$ (the $k$-subsets of $[t]$).
For two $k$-sets $X$ and $Y$ we set $W_i(X,Y)=1$ if and only if $|X \cap Y| = k-i$. Notice that $W_0=I$ and also notice that
$$
\sum_{i=0}^k W_i = J\;,
$$
where here $J$ denotes the all-$1$ matrix.
The set of matrices $J(t,k)$ is also known as the {\em Johnson Association Scheme}. The matrices of the Johnson scheme, as well as the
algebra formed by their linear
combinations, have been extensively studied. We refer the reader to \cite{VW} for an introduction to association
schemes, and the Johnson scheme in particular. For our purposes, we shall need the following explicit formulas for the eigenvalues of the matrices $W_i$ and their multiplicities.
\begin{lemma}
\label{l-eigen}
For each $0 \leq i \leq k$, the matrix $W_i$ has $k+1$ eigenvalues, denoted $p_i(0),p_i(1),\ldots,p_i(k)$.
\begin{itemize}
\item The multiplicity of $p_i(j)$ is ${t \choose j} - {t \choose {j-1}}$.
\item $p_i(j) = \sum_{r=0}^i (-1)^{i-r} {{k-r} \choose {i-r}}{{t-k+r-j} \choose r}{{k-j} \choose r}$.
\end{itemize}
\end{lemma}
Notice that, indeed, the sum of the multiplicities is ${t \choose k}$ and that $p_0(j)=1$ for all $j=0,\ldots,k$.

Another important property that we need is that any pair of matrices of $J(t,k)$ commute.
A classical result in linear algebra (see, e.g. \cite{HJ} Theorem 1.3.19) assets that if a set ${\cal S}$ of diagonalizable matrices
has the property that any pair of them commutes, then there exists a matrix $S$ which simultaneously diagonalizes each of them.
Namely, $SAS^{-1}$ is a diagonal matrix for any $A \in {\cal S}$.
Notice, that, in particular, this means that $S$ diagonalizes any linear combination of elements of ${\cal S}$.
In particular, we state another well known property of the Johnson Scheme:
\begin{lemma}
\label{l-comb}
If $C=\sum_{i=0}^k\alpha_iW_i$  then $C$ has $k+1$ eigenvalues $\lambda(0),\ldots,\lambda(k)$, where
for every $0 \leq j \leq k$ we have $\lambda(j) = \sum_{i=0}^k \alpha_i p_i(j)$, and the multiplicity of
$\lambda(j)$ is ${t \choose j} - {t \choose {j-1}}$.
\end{lemma}
We are now ready to prove Lemma \ref{hardpart}.

\paragraph{Proof of Lemma \ref{hardpart}:}
Let $C=A_{t,k,v}^TA_{t,k,v}$. Since, over the reals, $rank(X)=rank(X^TX)$, it suffices to prove that $rank(C) \geq  {t \choose k} - t+1$.
It is not difficult to see that $C$ is a linear combination of the elements of $J(t,k)$. Indeed, the rows and columns of
$C$ are indexed by ${{[t]} \choose k}$, and for two $k$-sets $X$ and $Y$, the value of $C(X,Y)$ is determined by $X \cap Y$.
Thus, $C= \sum_{i=0}^k\alpha_iW_i$. Clearly $\alpha_i$ simply counts the number of balanced $k$-cuts of $[t]$ for which
two $k$-sets $X$ and $Y$ with $|X \cap Y| = k-i$ are both transversals of the cut.
In order to better explain the main idea and main difficulty, we first consider the case $k=3$ (which will, in fact, hold for all $t \ge 12$).
This is the first non-trivial case of the lemma since the case $k=2$ follows from an easy lemma we prove at the end of this section.

The values of $\alpha_0,\alpha_1,\alpha_2,\alpha_3$ in the case $k=3$ are easily computed to be:
\begin{eqnarray*}
\alpha_0 & = & \frac{(t-3)!}{(t/3-1)!^3}\\
\alpha_1 & = & \frac{(t-4)!}{(t/3-1)!^2(t/3-2)!}\\
\alpha_2 & = & 2 \cdot \frac{(t-5)!}{(t/3-1)!(t/3-2)!^2}\\
\alpha_3 & = & 6 \cdot \frac{(t-6)!}{(t/3-2)!^3}\;.
\end{eqnarray*}
The eigenvalues of the matrices of $J(t,3)$ and their multiplicities are computed from Lemma \ref{l-eigen} and are given in Table \ref{table-1}.
\begin{table}
$$
\begin{array}{|c||c|c|c|c|c|}
\hline
j & p_0 & p_1 & p_2 & p_3 & {\rm multiplicity} \\
\hline
0 & 1 & -3 + 3(t-2) & 3-6(t-2)+3{{t-1} \choose 2} & -1+3(t-2)-3{{t-1} \choose 2}+{t \choose 3} & 1 \\
1 & 1 & -3 + 2(t-3) & 3-4(t-3)+{{t-2} \choose 2} & -1+2(t-3)-{{t-2} \choose 2} & t-1 \\
2 & 1 & -3 + (t-4) & 3-2(t-4) & -1 + (t-4) & {t \choose 2} - t \\
3 & 1 & -3 & 3 & -1 & {t \choose 3} - {t \choose 2} \\
\hline
\end{array}
$$
\caption{The eigenvalues of the matrices of $J(t,3)$ and their multiplicities}
\label{table-1}
\end{table}
We can now explicitly compute the eigenvalues of $C$ which are $\lambda(j) = \sum_{i=0}^k \alpha_i p_i(j)$ for $j=0,1,2,3$.
Clearly, $\lambda(0) > 0$ as for $t \ge 12$, $p_i(0) > 0$ for each $i=0,1,2,3$ (see Table \ref{table-1}).
It will be slightly more convenient to compute $\lambda(j)$ for $j=1,2,3$ by normalizing the $\alpha_i$, setting
$$
\alpha_i^* = \alpha_i \cdot \frac{(t/3-1)!^3}{(t-6)!}\;.
$$
Hence $\alpha_0^*=(t-3)(t-4)(t-5)$, $\alpha_1^*=(t-4)(t-5)(t/3-1)$, $\alpha_2^*=2(t-5)(t/3-1)^2$, $\alpha_3^*=6(t/3-1)^3$.
Denote the normalized eigenvalues by $\lambda(j)^*$. We obtain:
\begin{eqnarray*}
\lambda(1)^* & = & 1 \cdot [(t-3)(t-4)(t-5)]\\
           & + & [2t-9] \cdot [(t-4)(t-5)(t/3-1)]\\
           & + & [{{t-2} \choose 2}-4t+15] \cdot [2(t-5)(t/3-1)^2]\\
           & + & [2t-7 - {{t-2} \choose 2}] \cdot [6(t/3-1)^3]  = 0\;.
\end{eqnarray*}
Similarly,
\begin{eqnarray*}
\lambda(2)^* & = & 1 \cdot [(t-3)(t-4)(t-5)]\\
           & + & [t-7] \cdot [(t-4)(t-5)(t/3-1)]\\
           & + & [11-2t] \cdot [2(t-5)(t/3-1)^2]\\
           & + & [t-5] \cdot [6(t/3-1)^3]  > 0\;.
\end{eqnarray*}
$$
\lambda(3)^* = (t-3)(t-4)(t-5) -3(t-4)(t-5)(t/3-1) + 3 \cdot 2(t-5)(t/3-1)^2 - 6(t/3-1)^3  > 0.
$$
It follows that $\lambda(1)$ is the unique eigenvalue of $C$ whose value is $0$, and since its multiplicity is $t-1$ we obtain that
$rank(C) = {t \choose 3}-t+1$.

Note that the normalized eigenvalues of $C$ are polynomials in $t$ of degree bounded by a function of $k$ (actually, as we show below this degree is at most $2k$).
As is evident from the above proof of the case $k=3$, for large values of $k$, both the expressions for $p_i(j)$ in Lemma \ref{l-eigen} as well as the exact expressions
of $\alpha_i$, become very complicated. Hence it quickly\footnote{Actually, we have found this approach to be infeasible already for $k=4$.} becomes infeasible to {\em precisely} compute
the polynomials representing $\lambda^*(0),\ldots,\lambda^*(k)$.
Instead, we will use an {\em asymptotic} approach, by which
we will show that if we consider $\lambda^*(j)$ as a polynomial in $t$, then for every $j \neq 1$ the leading coefficient of this polynomial is positive. This
will imply that for all large enough $t$, the eigenvalues $\lambda(0),\lambda(2),\ldots,\lambda(k)$ are positive.
As the multiplicity of $\lambda(1)$ is $t-1$, this implies that $rank(C) \geq {t \choose k}-t+1$ which is what we need to show
\footnote{Note that even if we showed that
the leading coefficient of the polynomial representing $\lambda^*(1)$ is 0, it would not imply that
$\lambda^*(1)=0$. Hence, the proof only gives a lower bound for the rank of $C$. This lower bound is later matched by Lemma \ref{l-rank3}.}.

So we fix $k$, and assume, wherever necessary, that $t$ is sufficiently large.
We start with computing the expressions $\alpha_i$. We have
$$
\alpha_i = i! \cdot \frac{(t-k-i)!}{(t/k-1)!^{k-i} (t/k-2)!^i} \qquad {\mbox{for $i=0,\ldots,k.$}}
$$
To see this, recall that $\alpha_i$ counts the number of balanced $k$-cuts of $[t]$ for which
two $k$-sets $X$ and $Y$ with $|X \cap Y| = k-i$ are both transversals of the cut.

Again, it will be more convenient to normalize the $\alpha_i$ with
$$
\alpha_i^* = \alpha_i \cdot \frac{(t/k-1)!^k}{(t-2k)!} = i! (t/k-1)^i \prod_{s=k+i}^{2k-1}(t-s)\;.
$$
Notice that $\alpha_i^*$ is a polynomial in $t$ with degree $t^k$. Its leading coefficient is $i!/k^i$. Thus,
\begin{equation}
\label{e-alphai}
\alpha_i^* = \frac{i!}{k^i}t^k + O(t^{k-1})\;.
\end{equation}
Likewise, we can express $p_i(j)$ as a polynomial in $t$. From Lemma \ref{l-eigen} we obtain, for all $i=0,\ldots,k$, that
\begin{equation}
\label{e-pij1}
p_i(j) = \frac{t^i}{i!}{{k-j} \choose i}+O(t^{i-1}) \qquad\qquad\qquad\qquad\qquad {\mbox{for $j=0,\ldots,k-i.$}}
\end{equation}
\begin{equation}
\label{e-pij2}
p_i(j) = \frac{t^{k-j}}{(k-j)!}(-1)^{i-k+j}{j \choose {i-k+j}}+O(t^{k-j-1}) \qquad {\mbox{for $j=k-i+1,\ldots,k.$}}
\end{equation}
We are now ready to analyze $\lambda(j)^*= \sum_{i=0}^k \alpha_i^* p_i(j)$.
From (\ref{e-alphai}), (\ref{e-pij1}), and (\ref{e-pij2}) we immediately obtain that
$\lambda(j)^* = O(t^{2k-j})$. It therefore remains to show that $\lambda(j)^* = \Theta(t^{2k-j})$ for $j \neq 1$, that is,
that the coefficient of $t^{2k-j}$ does not vanish in these cases. From (\ref{e-alphai}), (\ref{e-pij1}), and (\ref{e-pij2}), this coefficient is:
$$
\frac{(k-j)!}{k^{k-j}} \cdot \frac{1}{(k-j)!} + \sum_{i=k-j+1}^k \frac{i!}{k^i}\cdot(-1)^{i-k+j}{j \choose {i-k+j}} \cdot \frac{1}{(k-j)!}\;.
$$
Rewriting it, we need to show that
$$
\sum_{s=0}^j (-1)^s \frac{1}{k^{s+k-j}}{j \choose s} \frac{(s+k-j)!}{(k-j)!} \neq 0\;.
$$
Notice that we know that the l.h.s. is nonnegative since $C$ is a positive semidefinite matrix,
and since the l.h.s. is a positive fraction of the leading coefficient of $\lambda(j)$ which is an eigenvalue of $C$.
Thus, equivalently, we must show that for every $j \in \{0,2,\ldots,k\}$ we have
\begin{equation}
\label{e-final}
\sum_{s=0}^j (-1)^s {j \choose s} k^{j-s} \frac{(s+k-j)!}{(k-j)!} > 0\;.
\end{equation}

Note that when $j=0$ the above expression is equal to 1 (for all $k$). Observe also
that the above is identically 0 when $j=1$, but again, this does not mean that $\lambda^{*}(1)=0$.
To prove (\ref{e-final}) for other values we will use the following result proved by Eli Berger \cite{Be}.
His proof involves a clever counting argument.
\begin{lemma}
\label{berger}
For integers $2 \le j \le k$:
$$
P_j(k) = \sum_{s=0}^j (-1)^s {j \choose s} k^{j-s} \frac{(s+k-j)!}{(k-j)!} > 0\;.
$$
\end{lemma}

\paragraph{Proof:}
Consider functions $f$ from $Z_j$ to $Z_k$.
We say that $f$ is {\em good at $i$}, if there exists $t \in \{0,\ldots,j-2\}$ so that $|f^{-1}\{f(i),f(i)+1,\ldots,f(i)+t\}| > t+1$.
Otherwise, $f$ is {\em bad at $i$}.
We say that $f$ is {\em good} if it is good at $i$ for all $i=0,\ldots,j-1$.
We claim that $P_j(k)$ counts the good functions.
Once we establish this claim notice that we are done since good functions exist, as any constant function is good
(taking $t=0$ for all $i$).
In order to prove that $P_j(k)$ counts the good functions we proceed as follows.
Let $F$ denote all the $k^j$ functions from $Z_j$ to $Z_k$, let $F_g \subset F$ denote the good functions,
and let $B_i \subset F$ denote the functions that are bad at $i$.
Clearly, $F_g = F \setminus \cup_{i=0}^{j-1} B_i$.
More generally, for as subset $U \subset \{0,\ldots,j-1\}$ let $B_U$ denote the functions that are bad for
all $i \in U$. In particular, $B_i = B_{\{i\}}$ and $B_\emptyset = F$. By the inclusion-exclusion principle,
$$
|F_g| = \sum_{U \subset \{0,\ldots,j-1\}} (-1)^{|U|}|B_U|\;.
$$
We will next prove that for $s=0,\ldots,j$, if $|U|=s$ then $|B_U| = k^{j-s}\frac{(k-j+s)!}{(k-j)!}$.
Once we prove this fact we have, by the last equation, that $|F_g|=P_j(k)$, as required.

First notice that for $s=0$ we trivially have $|B_\emptyset|= |F| = k^j$.
Now consider singletons $S=\{i\}$. We claim that for each of the $k^{j-1}$ possible assignments of values to $f(x)$ for $x \neq i$,
we can assign precisely $j-1$ values to $f(i)$ so as to obtain a function $f$ that is good at $i$ (and hence there are
$k^{j-1}(k-j+1)$ functions that are bad at $i$, as required).
Observe that the number of options of being good at $i$ is just a function
of the multiset of $j-1$ values at the $j-1$ points other than $i$.
Call a number $y \in Z_k$ in a multiset $Y$ (of elements of $Z_k$) {\em dense} if for some $t$ the set
$\{y, \ldots, y+t\}$ (modulo $k$) has more than $t+1$ elements in the multiset.
For a multiset $Y$, there are $|Y|$ ways to choose a number $y$ so that $y$ is dense in $Y+\{y\}$.
Hence, in our case, there are $j-1$ options for defining $f(i)$ so as to obtain a function that is good at $i$.

More generally, for subsets $U=\{u_1,\ldots,u_s\}$ of cardinality $s$, we have that for each of the $k^{j-s}$ possible assignments of values to $f(x)$
for locations $x$ with $x \notin S$, we can assign $j-s$ values to $f(u_1)$ so that $f$ will be good at $u_1$, and hence
$k-j+s$ values to $f(u_1)$ so that $f$ will be bad at $u_1$. Given such an assignment, there are now $k-j+s-1$ possible
values assigned to $f(u_2)$ so that $f$ is also bad at $u_2$. Similarly, having assigned values to $f(u_1),\ldots,f(u_x)$ so that
$f$ is bad at $u_y$ for $y=1,\ldots,x$, we can assign $k-j+s-x$ possible values to $f(u_{x+1})$ so that $f$ is bad also at
$u_{x+1}$. Overall, $|B_U| = k^{j-s}(k-j+s)(k-j+s-1)\cdots(k-j+1)$, as required.
\qed

\bigskip

Having verified (\ref{e-final}) the proof of Lemma \ref{hardpart} is now complete. \qed

\bigskip

We end this section with the proof of the following lemma which obtains another part of the statement of Theorem \ref{RankA}.
Recall that in order to avoid degenerate cases we assume that each coordinate of $v$ is at least $k$.

\begin{lemma}\label{fullrank}
If $v \neq (t/k,\ldots,t/k)$ then $rank(A_{t,k,v})={t \choose k}$.
\end{lemma}

As we have mentioned above the proof of the above lemma will be via a direct reduction to
Gottlieb's Theorem which we now turn to discuss.
For integers $t > h \ge k \ge 2$, the {\em inclusion matrix} $B(t,h,k)$ is defined as follows:
The rows of $B(t,h,k)$ are indexed by $h$-element subsets of $[t]$, and the columns by the
$k$-element subsets of $[t]$. Entry $(i,j)$ of $B(t,h,k)$ is $1$ if
the $k$-element set, whose index is $j$, is contained in the
$h$-element set, whose index is $i$. Otherwise, this entry is $0$.
Notice that $B(t,h,k)$ is a square matrix if and only if $t=h+k$,
and that for $t > h+k$, the matrix $B(t,h,k)$ has more rows than columns.
Trivially, $rank(B(t,h,k)) \leq {t \choose k}$. However, Gottlieb
\cite{G} proved that in fact
\begin{theo}[Gottlieb \cite{G}]
\label{gottlieb} $rank(B(t,h,k)) = {t \choose k}$ for all $t \ge h+k$.
\end{theo}

\paragraph{Proof of Lemma \ref{fullrank}:}
Since $v$ is not constant, we may assume that $v_1 > v_k$.
Also notice that since we always assume that $v_k \ge k$ we also have $t \ge t-v_k+k$.
Consider the set-inclusion matrix $B(t,t-v_k,k)$. It has the same number of columns as $A_{t,k,v}$. We will prove that the rows of $A_{t,k,v}$
span the rows of $B(t,t-a_k,k)$, and hence by Theorem \ref{gottlieb} we shall be done.

Fix a subset $T \subset [t]$ with $|T|=t-v_k$, and let $S_k = [t] - T$. Notice that $|S_k|=v_k$.
Consider any subset $S_1 \subset T$ with $|S_1|=v_1$ (if $k=2$ then $S_1=T$).
Consider any row $u$ of $A_{t,k,v}$ that corresponds to a partition $(S_1,S_2,\ldots,S_k)$ where $|S_i|=v_i$.
For every subset $S' \subset S_1$ of size $v_1-v_k$ let $u_{S'}$ be the row vector in $A_{t,k,v}$ corresponding to the
partition $(S' \cup S_k,S_2,\ldots,S_1-S')$.
We now observe that $u - u_{S'}$ has $+1$ for any $k$-set that is a transversal of $(S',S_2,\ldots,S_k)$,
and $-1$ for any $k$-set that is a transversal of $(S',S_2,\ldots,S_{k-1},S_1-S')$.
Therefore, if we take $w = \sum_{S' \subset S_1} u-u_{S'}$, where the summation is over all subsets of $S_1$ of size $v_1-v_k$,
we get a vector where all the $k$-sets that are a transversal of $(S_1,S_2,\ldots,S_k)$ have value ${{v_1-1} \choose {v_k}}$,
and all $k$-sets that have two vertices in $S_1$ and one vertex in each of $S_2,\ldots,S_{k-1}$ have value
$-{{v_1-2} \choose {v_k-1}}$.
An appropriate linear combination of $u$ and $w$ yields a vector which has $1$ only for $k$-sets with
two vertices in $S_1$ and one vertex in each of $S_2,\ldots,S_{k-1}$.
In any case, $k$-sets containing an element of $S_k$ get $0$.

Doing the same procedure for all possible choices of $S_1,S_2,\ldots, S_{k-1}$ shows that we can obtain the vector
which is constant on all $k$-sets that have all of their elements in $T$.
Hence, the row of $B(t,t-v_k,k)$ corresponding to $T$ is in the space spanned by the rows of $A_{t,k,v}$, as required.
$\qed$

\bigskip

We end this section with the discussion of the special case of $k=2$. Observe that when $v=(t/k,\ldots,t/k)$ the matrix $A_{t,k,v}$ contains
as a sub-matrix, the matrix $A_{t,k,v'}$ where we define $v'=(t/k,\ldots,t/k,t/k-1)$. Hence we immediately get from Lemma \ref{fullrank} that
$rank(A_{t,k,v}) \geq {t-1 \choose k}$. For $k \geq 3$ this bound does not match the bound we obtained in Lemma \ref{hardpart} but since
${t-1 \choose 2}={t \choose 2}-t+1$ this simple bound does show that Lemma \ref{hardpart} has a very simple proof for the special case $k=2$.

\section{The Structure of Counter-Examples}\label{secstruct}

In Sections \ref{Secunbalanced} and \ref{seccounter} we have given a characterization of the cut-properties which force a graph to be quasi-random.
In this section we will consider the non-quasi random cut properties, that is, the properties ${\cal P}_{\alpha}$ with
$\alpha=(1/r,\ldots,1/r)$. Our main result in this section will be a proof of Theorems \ref{theo2-informal} and \ref{theo3-informal} which will supply a description
of the hypergraphs satisfying ${\cal P}^{*}_{\alpha}$, which is the non-discrete version of ${\cal P}_{\alpha}$.
We will then use these results in order to derive an approximate description of all graphs satisfying ${\cal P}_{\alpha}$, see Theorem \ref{cutnorm}.
We will consider only the case $r=k$ as the proof for $r > k$ is identical.
Since throughout this section we assume that $\alpha=(1/k,\ldots,1/k)$ we will simplify the notation
by denoting ${\cal P}_{\alpha}$ and ${\cal P}^{*}_{\alpha}$ as ${\cal P}$ and ${\cal P}^{*}$, respectively.

Let us start by extending the definition of property ${\cal P}^{*}$, which was given in Subsection \ref{subsecstruct} for the special case of graphs, to the more general setting of hypergraphs.
In this case we are trying to assign weights to the edges of the complete $k$-uniform hypergraph on $t$ vertices, such that for any partition of the vertices into $k$ sets of sizes $t/k$ each, the
total weight of edges with exactly one vertex in each part is $p(t/k)^k$. Note that this problem can be cast in a linear algebra setting by trying to solve the following set
of linear equations. We have an unknown $x_s$ for every $k$-set of vertices $s$, and a linear equation $\ell_P$ for every partition of the vertices into $k$ sets of equal size.
If $P=(P_1,\ldots,P_k)$ then equation $\ell_P$ is $\sum_s x_s=p(t/k)^k$, where the sum is over all $k$-sets $s$ with exactly one vertex in each of the sets $P_1,\ldots,P_k$.
A key observation is that we can write this set of linear equations as $Ax=p(t/k)^k$ then $A$ is precisely the matrix $A_{t,k,v}$ defined in the previous section, where
$v=(t/k,\ldots,t/k)$.

So our goal now is to show that (assuming $t$ is large enough) the only solutions to the set of linear equations $A_{t,x,v}x=p(t/k)^k$ is the affine subspace spanned by the vectors corresponding
to $C_k(t,p)$ and $G_k(t,p)$ which were introduced in Subsection \ref{subseccharac}. More precisely, let $u_{t,k,p}$ be the ${t \choose k}$-dimensional vector all of whose entries
are $p$. This is the vector corresponding to the random $k$-uniform hypergraph $G_k(t,p)$ which satisfies ${\cal P}^{*}$. For a partition of $[t]$ into two sets of equal size $A,B$, let $v_{t,k,p}(A,B)$ be the
following ${t \choose k}$ dimensional vector; we think of the coordinates of $v_{t,k,p}(A,B)$ as being indexed by subsets of $[t]$ of size $k$. With this indexing we assign the entry of $v_{t,k,p}(A,B)$
corresponding the the set $\{v_{i_1},\ldots,v_{i_k}\}$ the value $2pj/k$ where $j=|\{v_{i_1},\ldots,v_{i_k}\} \cap A|$. Observe that this is exactly the vector representation of the hypergraphs
$C_k(t,p)$. As we have previously mentioned, we can actually define ${t-1 \choose t/2-1}$ such vectors, corresponding the possible ways of choosing the partition, and so we define
$V_{t,k,p}$ to be the collection of vectors $v_{t,k,p}(A,B)$ over all choices of $A,B$.

\begin{lemma}
\label{l-rank3}
The affine subspace spanned by the vector $u_{t,p,k}$ and $V_{t,p,k}$ has affine dimension at least $t$.
In particular, this implies that when $v=(t/k,\ldots,t/k)$ we have $rank(A_{t,k,v}) \leq {t \choose k}-t+1$.
\end{lemma}

\paragraph{Proof:}
Recall that $t$ vectors $w_1,\ldots,w_t$ are affine independent if and only if
the $t-1$ vectors $w_1-w_2,w_1-w_3, \ldots, w_1-w_t$ are linearly independent.
In our case, $w_1$ will represent the trivial solution (namely, the all-$p$ vector).

So, let $M$ be the matrix defined as follows. Its rows are all the ordered partitions of $[t]$ into two equal parts.
Its columns are all the $k$-tuples of $[t]$. For a row $(A,B)$ and a column $K$, the corresponding entry is $p-\frac{2p}{k}|B \cap K|$.
We need to prove that $rank(M) \ge t-1$. To simplify notation a bit we will divide each element of $M$ by $p$.
Hence, For a row $(A,B)$ and a column $K$, the corresponding entry is $1-\frac{2}{k}|B \cap K|$.

Let us first prove that a relatively small submatrix of $M$ already has a rank of $t-k+1$.
Let $U=\{1,\ldots,k-1\}$. Consider first the sub-matrix $M'$ of $M$ consisting only of rows $(A,B)$ where $U \subset A$,
and only of columns $K$ with $U \subset K$.
Notice that there are precisely $t-k+1$ such columns.
Notice that if we take each element $x$ of $M'$ and replace it with $(k/2)(x-1)$ we obtain
the inclusion matrix of singletons inside subsets of size $t/2$ of a $t-k+1$ element set.
In other words, we obtain the set-inclusion matrix $B(t-k+1,t/2,1)$.
Notice that $B(t-k+1,t/2,1)$ and $M'$ have the same rank since the sum of the of all the rows of $M'$ is a non-zero constant vector,
and hence the transformation $x \rightarrow (k/2)(x-1)$ on the elements of $M'$ does not change the rank.
By Theorem \ref{gottlieb}, $rank(B(t-k+1,t/2,1))=t-k+1$.
In particular, this means that there are $t-k+1$ rows of $M'$ that span $M'$ and form a basis to its rows.
Let $Z$ denote such $t-k+1$ rows, but now we think of them as rows of $M$ (not only $M'$).

Since $Z$ are independent restricted to the columns of $M'$, they are also independent as rows of $M$.
It remains to complement $Z$ with additional $k-2$ rows so as to form $t-1$ independent rows.

Consider the following $k-1$ rows, which we denote by $w_1,\ldots,w_{k-1}$.
Row $w_i$ corresponds to the partition $(A,B)$ where $A=\{i,t/2+2,\ldots,t\}$.
Notice that all the rows $w_1,\ldots,w_{k-1}$ are identical when restricted to the columns $K$ with $U \subset K$.
On the other hand, they are certainly not identical on the other columns. In fact, they are all independent.
Indeed, consider column $C_i$ where $C_i=\{i,t-k+2,\ldots,t\}$. Only $w_i$ has $1$ in this column and the other $w_j$ have $1-2/k$ in this column.
Thus, the $(k-1) \times (k-1)$ sub-matrix corresponding to the columns $C_i$ and to the rows $w_i$ is the all-$1$ matrix in the diagonal,
and $1-2/k$ anywhere else. This matrix is, of course, non-singular (it spans a non-zero constant vector and hence it also spans the same rows as the identity matrix). As the $w_i$ are identical when restricted to the columns $K$ with $U \subset K$, we see that by subtracting $w_{k-1}$ from each of the other $w_i$ we get equivalently that they span $k-2$ vectors $u_1,\ldots,u_{k-2}$ that are zero on the columns $K$ with $U \subset K$, and on the columns corresponding to $C_1,\ldots,C_{k-2}$ they form an $(k-2) \times (k-2)$ non-singular matrix.

We now get that $Z$, together with $u_1,\ldots,u_{k-2}$ form a set of $t-1$ independent vectors, that are all in the row space of $M$.
It follows that the row space of $M$ has rank at least $t-1$, as required.
\qed

\bigskip

The proof of Theorems \ref{theo2-informal} and \ref{theo3-informal} will now follow easily from the above result. We only prove Theorem \ref{theo3-informal} since it is clearly more general than Theorem \ref{theo2-informal}.

\paragraph{Proof of Theorem \ref{theo3-informal}:}
By Lemma \ref{hardpart}, when $v=(t/k,\ldots,t/k)$ we have $rank(A_{t,k,v}) \geq {t \choose k}-t+1$.
Hence, the affine subspace containing the solutions to $A_{t,k,v}\cdot x=p(t/k)^k$ has affine dimension at most $t$. But by Lemma
\ref{l-rank3} we get that the affine subspace spanned by the solutions corresponding to $G_k(t,p)$ and $C_k(t,p)$ has dimension at least
$t$. Hence, every solution to $A_{t,k,v}\cdot x=p(t/k)^k$ belongs to this subspace. $\qed$

\bigskip

We will now turn to show how to apply Theorem \ref{theo3-informal} regarding the solutions to ${\cal P}^{*}$ in order to obtain
Theorem \ref{cutnorm} below, which gives an approximate description of the hypergraphs satisfying ${\cal P}$.
Throughout the following we will ignore rounding issues as these have no effect on the asymptotic results.
In order to describe our result we need a few definitions.
We say that a $k$-uniform hypergraph $H=(V,E)$ with $n$ vertices is $\delta$-close to satisfying ${\cal P}$ if for every
partition of $V$ into $k$-sets of equal size $n/k$, the number of edges crossing the cut is $p(n/k)^k\pm \delta n^k$.
For a partition $P$ of $V$ into $t$ equal parts $V_1,\ldots,V_t$, we let $x=x_P$ denote the density-vector of $P$.
That is, $x$ has ${t \choose k}$ coordinates, indexed by the $k$-subsets of $[t]$, and $x_{K}$ is the density of the edges
of the cut induced by $(V_{i_1},\ldots,V_{i_k})$ where $K=\{i_1,\ldots,i_k\}$.
For integers $t$ and $k$, let $u_{t,p,k}$ and $V_{t,p,k}$ be the vectors that were defined before the statement of Lemma \ref{l-rank3}.
Recall that these vector encode the densities of $G_k(t,p)$ and $C_k(t,p)$.
In what follows we assume that $p$ and $k$ are constants.

\begin{theo}
\label{t-characterize}
For every $\epsilon > 0$ and $k \geq 2$, and for every large enough $t \geq t_0(k)$, there exists
$\delta=\delta(t,\epsilon) > 0$ so that the following holds for any $k$-uniform hypergraph. If $H$ is $\delta$-close to satisfying ${\cal P}$,
then for any partition $P$ of $V(H)$ into $t$ equal parts, the density vector $x_P$ satisfies $\ell_\infty(x_P,y) \leq \epsilon$ where
$y$ is an affine combination of $u_{t,p,k}$ and $V_{t,p,k}$.
\end{theo}

Note that the vector $y$ in the above theorem encodes a hypergraph on $t$ vertices satisfying ${\cal P}^{*}$. Thus the above theorem says that if we take any partition of the vertices of a hypergraph which is close to satisfying ${\cal P}$, then the densities of this partition are very close to the densities of a hypergraph satisfying ${\cal P}^{*}$.

\paragraph{Proof of Theorem \ref{t-characterize}:}
Suppose $H=(V,E)$ is an $n$-vertex hypergraph which is $\delta$-close to satisfying ${\cal P}$.
Fix any partition $P$ of $V$ into $t$ equal parts, $V_1,\ldots,V_t$.
Each partition ${\cal T}=\{Q_1,\ldots,Q_k\}$ of $[t]$ into $k$ equal parts corresponds to a partition of
$V_1,\ldots,V_t$ into $k$ equal parts $U_1,\ldots,U_k$, where $U_i = \cup_{j \in Q_i} V_j$.

As $U_1,\ldots,U_k$ is a balanced $k$-cut of $H$, we have that the number of edges of $H$ crossing this cut, denoted by $e_{{\cal T}}$,
satisfies
$$
\left| e_{{\cal T}} - p\left(\frac{n}{k}\right)^k \right| \le \delta n^k\;.
$$
It will be more convenient to write
\begin{equation}\label{eq1}
e_{{\cal T}} = p\left(\frac{n}{k}\right)^k + \rho_{{\cal T}}n^k,
\end{equation}
where $|\rho_{{\cal T}}| \le \delta$.

Another way to express $e_{{\cal T}}$ is via edge densities. For a $k$-subset $K \subset [t]$, let $d_K$ denote the density of edges
having one point in each $V_i$ where $i \in K$. Notice that the values $d_K$ are the entries of $x_P$.
Now, if ${\cal K}_{\cal T}$ is the set of $k$-subsets that are transversals of ${\cal T}$
(namely, have one point in each $Q_i$) then
\begin{equation}\label{eq2}
\sum_{K \in {\cal K}_{\cal T}} d_K\left(\frac{n}{t}\right)^k = e_{{\cal T}}\;.
\end{equation}
Combining (\ref{eq1}) and (\ref{eq2}) we get that for every ${\cal T}$ we have
$$
\sum_{K \in {\cal K}_{\cal T}} d_K = p\left(\frac{t}{k}\right)^k + t^k\rho_{{\cal T}}\;.
$$
Let $v=(t/k,\ldots,t/k)$ and suppose $t$ is large enough so that Theorem \ref{theo3-informal} holds, that is, that any solution to
$A_{t,k,v}\cdot x=p(t/k)^k J$ is an affine combination of $u_{t,p,k}$ and $V_{t,p,k}$ (here $J$ is the all-one vector).
Let $J'$ be a column vector, indexed by all the partitions ${\cal T}$ of $[t]$, where the entry of $J'$ corresponding to
${\cal T}$ is $t^k\rho_{{\cal T}}$. Observe that each partition of ${\cal T}$ is a row-index of $A_{t,k,v}$
and that each $K \subset [t]$ is a column index of $A_{t,k,v}$. It follows that $x_P$ is a solution
of the system
$$
A_{t,k,v}\cdot x =p\left(\frac{t}{k}\right)^kJ +J'\;.
$$
We already assume that $t$ is large enough so that each solution of $A_{t,k,v}\cdot x=p(t/k)^kJ$ is an affine combination of $u_{t,p,k}$ and $V_{t,p,k}$.
Each element of $J'$ has absolute value at most $t^k\delta$, hence $||J'||_{\infty}$ converges to zero with $\delta$.
It is easy to see (as we show in the next paragraph) that as $||J'||_{\infty}$ converges to zero, any solution to
$A_{t,k,v} \cdot x =\left(\frac{t}{k}\right)^kJ+J'$ converges to a solution
of $A_{t,k,v} \cdot x =\left(\frac{t}{k}\right)^kJ$. In particular, it follows that for $\delta$ sufficiently small, $x_P$ is $\epsilon$-close to an
affine combination of $u_{t,p,k}$ and $V_{t,p,k}$, and the result follows.

For completeness, we show that for every $\epsilon' > 0$ and for every matrix $A$, there is $\delta'=\delta'(A,\epsilon') > 0$ so that if
$||b-b'||_{\infty} < \delta'$ then for any solution $x_1$ of the system $Ax=b'$ there exists a solution $x_2$ of the system $Ax=b$ so that
$||x_1-x_2||_{\infty} < \epsilon'$. Let $d=rank(A)$. Notice that we may assume that $A$ has full row rank (namely $A$ is some $d \times n$ matrix),
since we may always truncate ``unnecessary'' rows (of $A$ and of $b$ and $b'$) and the solution spaces of the truncated systems remain intact.
We may now assume, without loss of generality that $A=[A_1|A_2]$ where $A_1$ is a non-singular $d \times d$ matrix consisting of the first $d$ columns,
and $A_2$ is the remaining $n-d$ columns. So let $x_1$ be a solution to $Ax=b'$. Let $y_1$ be the truncation of $x_1$ to the first $d$ entries,
and let $z_1$ be the truncation of $x_1$ to the last $n-d$ entries. We now have that $y_1$ is the {\em unique} solution to the system
$A_1x=b'-A_2z_1$. In other words, $y_1$ is just $A_1^{-1}(b'-A_2z_1)$. Now let $y_2$ be the unique solution of $A_1x=b-A_2z_1$.
Notice that since $||b-b'||_{\infty} < \delta'$ then trivially also $||(b-A_2z_1)-(b'-A_2z_1)||_{\infty} < \delta'$. Since the mapping $v \rightarrow A_1^{-1}v$ is continuous,
we have that for $\delta'$ sufficiently small, $||y_1-y_2||_{\infty} < \epsilon'$. But now define $x_2$ to be the vector whose first $d$ coordinates are $y_2$
and whose last $n-d$ coordinates are $z_1$. Notice that $Ax_2=b$ and $||x_1-x_2||_{\infty}=||y_1-y_2||_{\infty} < \epsilon'$.
\qed

\bigskip

Let $G,G'$ be two (possibly weighted) graphs. A natural and well studied measure for the distance between two graphs is the
cut-norm introduced by Frieze and Kannan \cite{FK}. We briefly mention that the cut-norm is central to the study of graph limits and refer the reader to \cite{LSz} for more information and references.
Let's start with the basic definitions. We denote by $e_{G}(S,T)$ the total weight of the
edges of $G$ connecting $S$ to $T$, where edges belonging to both $S$ and
$T$ are counted twice. The cut-norm between two graphs $G$ and $G'$ on a set of $n$ vertices is then defined to be
$$
d_{\Box}(G,G')=\frac{1}{n^2}\max_{S,T \subseteq [n]}|e_G(S,T)-e_{G'}(S,T)|\;.
$$
The following is our approximate description of the graphs satisfying ${\cal P}$.

\begin{theo}\label{cutnorm}
for every $\epsilon > 0$ there is a $\delta=\delta(\epsilon)>0$ such that if $G$ is $\delta$-close to satisfying ${\cal P}$, then there is a graph
$G'$ satisfying ${\cal P}^{*}$ for which $d_{\Box}(G,G')\leq \epsilon$. Moreover,
$G'$ has {\em constant complexity}; it is an affine combination of $G(n,p)$ and a number of copies of $C_2(n,p)$ which depends
only on $\epsilon$.
\end{theo}

For the proof of Theorem \ref{cutnorm} we will need the so called {\em weak regularity lemma} of Frieze and Kannan \cite{FK}. To state
this lemma we need the following notation. An {\em equipartition} $P=\{V_1,\ldots,V_t\}$ of a graph $G=(V,E)$ is a partition of $V$ into
subsets of equal size. The {\em order} of the equipartition is the number of sets in $P$. Given a graph $G$ and an equipartition
$P$ we define $G[P]$ to be the following weighted graph. If vertices $u$ and $v$ both belong to one of the sets
of $P$ then the weight of the edge $(u,v)$ is set to zero. Otherwise, there are $i<j$ such that $v \in V_i$ and $u \in V_j$
and in this case we assign the edge $(u,v)$ a weight $e(V_i,V_j)/|V_i||V_j|$, where $e(V_i,V_j)$ denotes the number of edges connecting
$V_i$ and $V_j$. The result of Frieze and Kannan \cite{FK} can be stated as follows.

\begin{theo}[Frieze and Kannan \cite{FK}]\label{FKtheo} For every $\epsilon > 0$ there is an integer $T=T(\epsilon)$ satisfying the following. Every graph $G$
has an equipartition $P$ of order $1/\epsilon \leq t \leq T$ satisfying $d_{\Box}(G,G[P]) \leq \epsilon$.
\end{theo}

We note that the above theorem can also be deduced from the regularity lemma of Szemer\'edi \cite{Sz}. However, while the bound on $T(\epsilon)$
in Theorem \ref{FKtheo} grows like $2^{O(1/\epsilon^2)}$, the bound one obtains from Szemer\'edi's regularity lemma are significantly weaker.

\paragraph{Proof of Theorem \ref{cutnorm}:} Given $\epsilon > 0$ let $T=T(\epsilon/2)$ be the constant from Theorem \ref{FKtheo}
and set $\delta=\delta(T,\epsilon/4)$ to be the constant from Theorem \ref{t-characterize}. Suppose $G=(V,E)$ is $\delta$-close to satisfying ${\cal P}$. Denote $V(G)$ by $[n]$. Applying Theorem \ref{FKtheo} on $G$ with $\epsilon/2$ we obtain an equipartition of $[n]$ of order $2/\epsilon \leq t \leq T$ satisfying $d_{\Box}(G,G[P]) \leq \epsilon/2$. Since $d_{\Box}$ satisfies the triangle inequality we finish the proof by showing that there is a graph $G'$ satisfying ${\cal P}^*$ and $d(G[P],G') \leq \epsilon/2$. Let $x_P$ be the density vector corresponding to $P$. Since we assume
that $G$ is $\delta$-close to satisfying ${\cal P}$ we get from the choice of $\delta$ and Theorem \ref{t-characterize}, that there exists
a density vector $y$ satisfying $\ell_{\infty}(x_P,y) \leq \epsilon/4$, where $y$ is an affine combination of the vectors $u_{t,p,2}$ and $V_{t,p,2}$.
We now wish to use $y$, which encodes a graph on $t$ vertices satisfying ${\cal P}^{*}$, in order to define a graph
on $[n]$ satisfying ${\cal P}^{*}$. We do this in the most ``obvious'' way.
First, we can replace $u_{t,p,2}$ which encodes a random graph on $t$ vertices, with a random graph on $n$ vertices. As to the vectors of $V_{t,p,2}$, for every $v \in V_{t,p,2}$ we define a weighted graph $G_v$ on $V(G)$ as follows; if $u$ and $w$ belong to one of the sets
$V_1,\ldots,V_t$ we assign $(u,w)$ weight zero. Otherwise, there is $i<j$ such that $u \in V_i$ and $w \in V_j$ in which case we assign
$(u,w)$ the weight assigned by $v$ to the pair $(i,j)$. Observe that for every $v$ the graph $G_v$ is simply a copy of $C_2(n,p)$.
Now recall that we assume that $y$ is an affine combination of $u_{t,p,2}$ and the vectors of $V_{t,p,2}$. If $y$ can be written as $y=\lambda\cdot u_{t,p,2}+\sum_{v \in V_{t,p,2}}\lambda_{v}\cdot v$ then we set $G'=\lambda \cdot G(n,p)+\sum_{v \in V_{t,p,2}}\lambda_{v}\cdot G_v$. Observe that by defining $G'$ this way we guarantee that for every pair of vertices $u \in V_i$ and $w \in V_j$ the weight assigned to $(u,w)$ in $G'$ is the weight
assigned to $(i,j)$ by $y$. Furthermore, by the properties of $y$, we know that this is within $\epsilon/4$ of the weight assigned to $(u,w)$ in
$G[P]$.

To see that $d(G[P],G') \leq \epsilon/2$ consider any pairs of sets $S,T \subseteq [n]$. Since for every $i < j$ and every $u \in V_i$ and $w \in V_j$ the weights assigned to the edge $(u,w)$ in $G'$ and $G[P]$ differ by at most $\epsilon/4$ we get that these edges contribute to $|e_{G[P]}(S,T)-e_{G'}(S,T)|$
at most $\frac14\epsilon|S||T| \leq \frac14\epsilon n^2$. The other contribution is due to edges which belong to one of the sets $V_i$. But since $t \geq 2/\epsilon$ the total
contribution of such edges is at most $t {n/t \choose 2} \leq \frac14\epsilon n^2$. All together we get that $d(G[P],G') \leq \epsilon/2$
as needed. Finally, note the graph $G'$ is indeed obtained as
an affine combination of a number of copies of $C_2(n,p)$ which depends only on $\epsilon$ and is independent of the size of $G$. $\qed$

\section{An Open Problem}\label{secopen}

In this paper we extended the Chung-Graham result \cite{CG} from $2$-cuts of graphs to arbitrary cuts of graphs and hypergraphs.
We would like to raise as an open problem the possibility of obtaining another extension of the Chung-Graham result. Motivated
by the results of Simonovits and S\'os \cite{SSnind,SSind}, there have been several recent investigations \cite{CHPS,DR,S,SY,Y} which suggest that in some sense, many of the
quasi-random properties regarding edge distributions, remain quasi-random if one replaces an edge by any other fixed graph. For example,
property ${\cal P}_1$, which defines what it means for a graph to be quasi-random, requires the number of edges to be the ``correct'' one in all subsets of vertices. A theorem of Simonovits and S\'os \cite{SSnind} asserts that ${\cal P}_1$ is actually equivalent to property ${\cal P}_H$ in which the edge is replaced by a fixed graph $H$. Another example, is property ${\cal P}_2$ in Theorem \ref{CGtheo} which asserts that in order to guarantee that $G$
is quasi-random it is enough to require that only sets of size $\alpha n$ have the correct number of edges. A similar variant of
the Simonovits and S\'os \cite{SSnind} theorem mentioned above was obtained recently in \cite{S,Y}. Given the above discussion, it is natural to consider the number of copies of a fixed graph $H$ that have one vertex in each of the classes of a cut. For simplicity, let's first consider cliques
of size $k$ and $k$-cuts\footnote{Actually, as we show
in the main body of the paper, the case of cuts with more classes can be reduced to the special case where the number of classes equals the number
of vertices of $K_k$.}.

\begin{definition}[$C_{\alpha}$] Let $\alpha=(\alpha_1,\ldots,\alpha_k)$ be a vector of positive reals satisfying $\sum_i \alpha_i=1$. Given a partition of $V(G)$ into $k$ sets of sizes $|V_i|=\alpha_i n$ we denote by $C(V_1,\ldots,V_k)$ the number of copies of $K_k$ in $G$ with precisely one vertex
in each of the sets $V_i$. We say that a graph satisfies $C_{\alpha}$ if for any partition of the vertices of $V(G)$ into $k$ sets of sizes $|V_i|=\alpha_i n$ we have
\begin{equation}\label{defC}
C(V_1,\ldots,V_k)=\frac{p^{{k \choose 2}}}{k!}n^k\prod^k_{i=1}\alpha_i \pm o(n^k)\;.
\end{equation}
\end{definition}

So we now ask which cut properties ${\cal C}_{\alpha}$ are quasi-random? The following is a simple corollary of our main result and
the theorem of Simonovits and S\'os \cite{SSnind} mentioned above.

\begin{prop}\label{easyprop} If $\alpha \neq (1/k,\ldots,1/k)$ then $C_{\alpha}$ is quasi-random.
\end{prop}

\paragraph{Proof:} Suppose $G$ satisfies $C_{\alpha}$ and define a $k$-uniform hypergraph $H$ which has an edge on a $k$-tuple of vertices $v_1,\ldots,v_k$ if and only if
$v_1,\ldots,v_k$ form a clique in $G$. The assumption that $G$ satisfies (\ref{defC}) in any $\alpha$-cut
means that $H$ satisfies the $\alpha$-cut property ${\cal P}_{\alpha}$ of $k$-uniform hypergraphs with edge density $p^{{k \choose 2}}$. Since
$\alpha \neq (1/k,\ldots,1/k)$ we get from Theorem \ref{theomain} that $H$ is $p^{{k \choose 2}}$-quasi-random. Going back to the graph $G$, this means that every set of vertices $U$ in $G$ has the correct number of copies of $K_k$ we expect to find in $G(n,p)$. Therefore, by the result of Simonovits and S\'os \cite{SSnind} mentioned above this means that $G$ is $p$-quasi-random, thus completing the proof. $\qed$

\bigskip

So the above proposition shows that for non-balanced $\alpha$-cuts, the property ${\cal C}_{\alpha}$ is quasi-random via a reduction to property ${\cal P}_{\alpha}$ of Theorem
\ref{theomain}. Since ${\cal P}_{\alpha}$ is not quasi-random for balanced $\alpha$-cuts, this still leaves the case of balanced cuts open. We thus raise the following open problem.

\begin{problem}\label{openproblem} Is the property $C_{\alpha}$ quasi-random when $\alpha=(1/k,\ldots,1/k)$.
\end{problem}

One can of course wonder why is it the case that the example from Section \ref{seccounter} showing that ${\cal P}_{\alpha}$ is not quasi-random for balanced $\alpha$, does not imply that $C_{\alpha}$ is also not quasi-random. The reason is that there is no obvious way of defining a graph, whose copies of $K_k$ will have the same distribution as the edges of the hypergraph that shows that ${\cal P}_{\alpha}$ is not quasi-random.
Actually, we can prove that one cannot construct a counter-example to Problem \ref{openproblem} by ``imitating'' the construction we used to show that ${\cal P}_{\alpha}$ is not quasi-random. More precisely, recall
that for $k=3$ the counter example, denoted $C_3(n,p)$, was obtained by partitioning the vertices into two sets $A,B$ of equal size and putting an edge containing the set of vertices $\{v_1,v_2,v_3\}$ with probability proportional to $|\{v_1,v_2,v_3\} \cap A|$. It is thus natural to ask if one can find three numbers $p_1,p_2,p_3$ such that if one picks every edge in $A$ with probability $p_1$, every edge in $B$ with probability $p_2$ and every edge between $A$ and $B$ with probability $p_3$, then the distribution of triangles in this graph will be identical to the distribution of edges in $C_3(n,p)$. It is not hard to to see that by applying Theorem \ref{gottlieb}
one can prove that such $p_1,p_2,p_3$ do not exist.

Recall that the results described in Section \ref{secstruct} give a description of all hypergraphs satisfying ${\cal P}_{\alpha}$.
From this result it follows that if one can construct a counter example showing that the answer to Problem \ref{openproblem} is negative then the distribution of triangles in this example would
have to be an affine combination of the distributions of $3$-edges in $C_3(n,p)$ and the random 3-uniform hypergraph with edge probability $p$. As we have just argued we can show that such an example cannot be obtained by imitating the distribution of edges in a ``single'' copy of $C_3(n,p)$.
However, this does not rule out the possibility of constructing such an example by imitating the distribution of edges in a {\em combination} of several copies of $C_3(n,p)$. It seems very interesting to further investigate this problem.

\paragraph{Acknowledgements:} We would like to thank Nati Linial for very helpful discussions and Eli Berger for his proof of Lemma \ref{berger}.


\begin{thebibliography}{99}

\bibitem{ADLRY}
N. Alon, R. A. Duke, H. Lefmann, V. R\"odl and R. Yuster, The
algorithmic aspects of the regularity lemma, J. Algorithms 16
(1994), 80-109.

\bibitem{BF}
L. Babai and P. Frankl, {\bf Linear Algebra Methods in Combinatorics}, book manuscript, 1992.

\bibitem{Be}
E. Berger, Private communication, 2009.

\bibitem{C}
F. R. K. Chung, Regularity lemmas for hypergraphs and
quasi-randomness, Random Structures and Algorithms 2 (1991),
241-252.

\bibitem{C2}
F. R. K. Chung, Quasi-random classes of hypergraphs, Random Structures Algorithms 1
(1990), 363–382.

\bibitem{CG}
F. R. K. Chung and R. L. Graham, Quasi-random set systems, Journal
of the AMS, 4 (1991), 151-196.

\bibitem{CG1}
F. R. K. Chung and R. L. Graham, Quasi-random tournaments, J.
Graph Theory 15 (1991), 173-198.


\bibitem{CG2}
F. R. K. Chung and R. L. Graham, Quasi-random hypergraphs, Random
Structures and Algorithms 1 (1990), 105-124.

\bibitem{CG3}
F. R. K. Chung and R. L. Graham, Maximum cuts and quasi-random
graphs, Random Graphs, (Poznan Conf., 1989) Wiley-Intersci, Publ.
vol 2, 23-33.

\bibitem{CGW}
F. R. K. Chung, R. L. Graham and R. M. Wilson, Quasi-random graphs,
Combinatorica 9 (1989), 345-362.

\bibitem{CHPS}
D. Conlon, H. H\`an, Y. Person and M. Schacht,
Weak quasi-randomness for uniform hypergraphs, submitted, 2009.

\bibitem{DR}
D. Dellamonica, Jr. and V. R\"odl, Hereditary quasi-random properties of hypergraphs, submitted.

\bibitem{FK}
A. Frieze and R. Kannan, Quick approximation to matrices and
applications, {\em Combinatorica} 19 (1999), 175-220.

\bibitem{FR}
P. Frankl and V. R\"odl, Extremal problems on set systems, Random
Structures and Algorithms 20 (2002), 131-164.


\bibitem{G}
D. H. Gottlieb, A class of incidence matrices, Proc. Amer. Math.
Soc. 17 (1966), 1233-1237.

\bibitem{Go}
T. Gowers, Quasirandom groups, Combinatorics, Probability and Computing 17 (2008), 363-387.

\bibitem{Go1}
T. Gowers, Quasirandomness, counting and regularity for 3-uniform
hypergraphs, Combinatorics, Probability and Computing 15 (2006), 143-184.

\bibitem{Go2}
T. Gowers, Hypergraph regularity and the multidimensional
Szemer\'edi theorem, Ann. of Math. 166 (2007), 897-946.

\bibitem{HLW}
S. Hoory, N. Linial and A. Wigderson, Expander graphs and their
applications, Bulletin of the AMS, Vol 43 (4), 2006, 439-561.

\bibitem{HJ}
R.A. Horn and C.R. Johnson,
Matrix Analysis,
Cambridge University Press, London, 1990.


\bibitem{Ish}
Y. Ishigami, A simple regularization of hypergraphs, at:
http://arxiv.org/abs/math/0612838.

\bibitem{J}
S. Janson, Quasi-random graphs and graph limits, manuscript, 2009.

\bibitem{KNRS}
Y. Kohayakawa, B. Nagle, V. R\"odl, and M. Schacht, Weak hypergraph regularity and linear
hypergraphs, J. Comb. Theory Ser. B, to appear.


\bibitem{KS}
M Krivelevich and B. Sudakov,  Pseudo-random graphs, More sets,
graphs and numbers, E. Gy\H{o}ri, G. O. H. Katona and L. Lov\'asz, Eds.,
Bolyai Society Mathematical Studies Vol. 15, 199-262.

\bibitem{LS}
L. Lov\'asz and V. T. S\'os, Generalized quasirandom graphs,
Journal of Combinatorial Theory Series B 98 (2008), 146-163.

\bibitem{LSz}
L. Lov\'asz and B. Szegedy, Szemer\'edi's Lemma for the analyst, GAFA 17 (2007), 252-270.

\bibitem{NRS}
B. Nagle, V. R\"odl and M. Schacht, The counting lemma for regular
$k$-uniform hypergraphs, Random Structures and Algorithms 28
(2006), 113-179.

\bibitem{NPRS}
B. Nagle, A. Poerschke, V. R\"odl, and M. Schacht, Hypergraph regularity and quasirandomness,
Proc. of the $20^{th}$ Annual ACM-SIAM Symposium on Discrete Algorithms (2009) 227–235.



\bibitem{RS}
V. R\"odl and J. Skokan, Regularity lemma for $k$-uniform
hypergraphs, Random Structures and Algorithms 25 (2004), 1-42.


\bibitem{S}
A. Shapira, Quasi-randomness and the distribution of copies of a
fixed graph, Combinatorica, 28 (2008), 735-745.

\bibitem{SY}
A. Shapira and R. Yuster, The effect of induced subgraphs on quasi-randomness, Proc. of SODA 2008, 789-798.
Also, Random Structures and Algorithms, to appear.

\bibitem{SS}
M. Simonovits and V. T. S\'os, Szemer\'edi's partition and
quasirandomness, Random structures and algorithms, 2 (1991), 1-10.

\bibitem{SSnind}
M. Simonovits and V. T. S\'os, Hereditarily extended properties,
quasi-random graphs and not necessarily induced subgraphs,
Combinatorica 17 (1997), 577-596.

\bibitem{SSind}
M. Simonovits and V. T. S\'os, Hereditarily extended properties,
quasi-random graphs and induced subgraphs, Combinatorics Probability
and Computing, 12 (2003), 319-344.


\bibitem{Sztheo}
E.~Szemer\'edi, Integer sets containing no $k$ elements in
arithmetic progression, Acta Arith. 27 (1975), 299-345.

\bibitem{Sz} E.~Szemer\'edi,
Regular partitions of graphs, In: {\em Proc.\ Colloque Inter.\
CNRS} (J.~C.~Bermond, J.~C.~Fournier, M.~Las~Vergnas and
D.~Sotteau, eds.), 1978, 399--401.

\bibitem{T}
T. Tao, A variant of the hypergraph removal lemma, J.
Combin. Theory, Ser. A 113 (2006), 1257-1280.

\bibitem{T1}
A. Thomason, Pseudo-random graphs, Proc. of Random Graphs, Pozna\'n
1985, M. Karo\'nski, ed., Annals of Discrete Math. 33 (North Holland
1987), 307-331.

\bibitem{T2}
A. Thomason, Random graphs, strongly regular graphs and
pseudo-random graphs, Surveys in Combinatorics, C. Whitehead, ed.,
LMS Lecture Note Series 123 (1987), 173-195.

\bibitem{VW}
J.H. van Lint and R.M. Wilson,
A course in Combinatorics, Cambridge University Press, London, 1992.

\bibitem{Y}
R. Yuster, Quasi-randomness is determined by the distribution of copies of a fixed graph in equicardinal large sets,
Proc. of APPROX-RANDOM 2008, 596-601.




\end{thebibliography}
\end{document}